\documentclass[12pt]{article}

\usepackage{amscd, amssymb, amsthm, amsmath,eucal, verbatim}

\newcommand{\proj}{\mathbf P}

\newcommand{\rarr}{\rightarrow}

\newcommand{\com}{\mathbb{C}}
\newcommand{\Q}{{\mathbb{Q}}}
\newcommand{\Z}{{\mathbb{Z}}}

\newcommand{\ff}{{\mathbf{f}}}
\newcommand{\gb}{{\mathbf{g}}}
\newcommand{\pp}{{\mathbf{p}}}
\newcommand{\ee}{{\mathbf{e}}}
\newcommand{\zz}{{\mathfrak{z}}}
\newcommand{\cs}{{\varsigma}}
\newcommand{\Th}{\vartheta} 
\newcommand{\Hr}{\mathsf{H}}

\newcommand{\wtko}{\widetilde{(k+1)}}
\newcommand{\wtk}{\widetilde{(k)}}
\newcommand{\wtp}{\widetilde{\pp}}

\newcommand{\la}{\lambda}

\newcommand{\cE}{\mathcal{E}}
\newcommand{\cP}{\mathcal{P}}
\newcommand{\cS}{\mathcal{S}}
\newcommand{\cF}{\mathcal{F}}

\newcommand{\al}{\alpha}
\newcommand{\vac}{v_\emptyset}

\newcommand{\lang}{\left\langle}
\newcommand{\rang}{\right\rangle}

\newcommand{\gli}{\mathfrak{gl}(\infty)}
\newcommand{\nr}[1]{:\!#1\!:}
\newcommand{\ul}{\underline}
\newcommand{\LV}{\Lambda^{\frac\infty2}V}
\newcommand{\LVc}{\Lambda^{\frac\infty2}_0V}
\newcommand{\Ls}{\Lambda^*}
\newcommand{\sh}{{\textstyle \frac12}}

\newcommand{\lam}{\lambda}

\newcommand{\ZZ}{\mathcal Z}

\DeclareMathOperator{\ev}{ev} \DeclareMathOperator{\Aut}{Aut}

\newcommand{\C}{\mathbb{C}}

\DeclareMathOperator{\tr}{tr}

\DeclareMathOperator{\gr}{gr}

\newtheorem{pr}{Proposition}[section]
\newtheorem{lm}[pr]{Lemma}
\newtheorem{tm}{Theorem}
\newtheorem{cor}[pr]{Corollary}

\theoremstyle{definition}

\numberwithin{equation}{section}

\begin{document}
\title{Gromov-Witten theory, Hurwitz theory,
and completed cycles.}
\author{A.~Okounkov
and R.~Pandharipande}
\date{April 2002}
\maketitle

\setcounter{tocdepth}{2}
\tableofcontents



\setcounter{section}{-1}
\section{Introduction}
\subsection{Overview}

\subsubsection{}

There are two enumerative theories
of maps from 
curves to curves.
Our goal here is to study their relationship.
All curves in the paper
will be projective over $\com$.

The first theory, introduced in the $19^{th}$ century by
Hurwitz, concerns the
enumeration of degree $d$ covers, 
$$\pi:C\to X,$$ 
of nonsingular curves $X$ with 
specified ramification data.  In 1902,
Hurwitz published a  closed formula for the number of
covers, $$\pi:\proj^1\rarr \proj^1,$$ with specified simple
ramification over ${\mathbf A}^1\subset \proj^1$ and arbitrary
ramification over $\infty$
(see \cite{Hur} and also  \cite{ELSV,Str}).
 
Cover enumeration is easily expressed in the
class algebra of the symmetric group $S(d)$. The formulas
involve 
the characters of $S(d)$.
Though great strides have been taken in the past century,
the characters of $S(d)$ remain 
objects of substantial combinatorial complexity. 
While any particular Hurwitz number may be calculated,
very few explicit formulas are available.

The second theory,
the Gromov-Witten theory of 
target curves $X$, is
modern.
Gromov-Witten theory  
is defined via intersection in the  
moduli space $\overline{M}_{g,n}(X,d)$
of degree $d$ stable maps, $$\pi:C\to X,$$ 
from genus $g$, $n$-pointed
curves. 
A sequence of {descendents}, 
\begin{equation*}
\tau_0(\gamma), \tau_1(\gamma), \tau_2(\gamma) ,\ldots\, ,
\end{equation*}
is determined
by each 
cohomology class $\gamma\in H^*(X,\Q)$. The descendents $\tau_k(\gamma)$ 
correspond to classes in the
cohomology
of $\overline{M}_{g,n}(X,d)$.
Full definitions are
given in Section \ref{secgw} below.  
The Gromov-Witten invariants of
$X$ are defined as integrals of products of
descendent classes against the
virtual fundamental class
of $\overline{M}_{g,n}(X,d)$.

Let 
$\omega\in H^2(X,\Q)$ denote the (Poincar\'e dual) class of a point.
We define the {\em stationary  sector} of the Gromov-Witten theory $X$
to be the integrals
involving only the descendents of $\omega$.
The stationary sector is the most basic and fundamental part of the
Gromov-Witten theory of $X$. 

Since Gromov-Witten theory  and  Hurwitz theory are both
enumerative theories of maps,
we may ask 
whether there is any precise relationship between the two. 
We prove
the stationary sector of Gromov-Witten
is in fact {\em equivalent} to Hurwitz theory.

\subsubsection{}

Let $X$ be a nonsingular target curve.
The main result of the paper is a 
correspondence, termed here the {\em GW/H correspondence}, 
between the stationary sector of  Gromov-Witten theory  and
Hurwitz theory.

Each 
descendent $\tau_k(\omega)$ corresponds to 
an explicit linear combination of  
ramification conditions in Hurwitz theory.
A stationary Gromov-Witten invariant of $X$ is equal to the sum 
of the Hurwitz numbers obtained
by replacing $\tau_k(\omega)$ by
the associated ramification conditions.
The ramification conditions associated to $\tau_k(\omega)$
are universal ---
independent of all factors including the target $X$.

\subsubsection{}

The GW/H correspondence
may be alternatively expressed as associating to each descendent
$\tau_k(\omega)$ an 
explicit element of the class algebra of the symmetric group. 
The associated elements, the completed cycles,
have been considered previously in
Hurwitz theory --- the term {\em completed cycle} first appears in
 \cite{EOZ} following unnamed appearances of the associated elements
in \cite{BO,EO}.
In fact, completed cycles, implicitly, 
are ubiquitous in the theory of shifted symmetric functions. 

The completed $k$-cycle
is the ordinary $k$-cycle 
corrected by a non-negative linear combination
of permutations with smaller support 
(except, possibly, for the constant term corresponding to the
empty permutation, which may be of either sign).
The corrections
are viewed as completing the cycle.
In \cite{EOZ},
the corrections to the ordinary $k$-cycle were understood 
as
counting degenerations of Hurwitz coverings with
appropriate combinatorial weights. Similarly,
in Gromov-Witten theory, the correction terms will be seen
to arise from the
boundary strata of $\overline{M}_{g,n}(X,d)$.

\subsubsection{}

The GW/H correspondence is important from several 
points of view. From the geometric perspective, the correspondence
provides a combinatorial approach to the stationary Gromov-Witten
invariants of $X$, leading to very concrete and efficient formulas.
 From the perspective of symmetric functions,
a geometrization of the theory of completed cycles is obtained.

Hurwitz 
theory with completed cycles is combinatorially 
much more accessible than  
standard Hurwitz theory --- a major motivation
for the 
introduction of completed cycles. 
Completed cycles calculations may be naturally evaluated in the
operator formalism of  
the infinite wedge representation, $\LV$.
In particular,
closed formulas for the completed cycle correction
terms are obtained.
If the target 
$X$ is either genus 0 or 1, closed
form evaluations of all corresponding generating functions
may be found, see Sections \ref{sP1} and \ref{sE}. 
In fact, the completed cycle corrections appear in the
theory with target genus 0.

Hurwitz theory, while elementary
to define, leads to substantial combinatorial difficulties.
Gromov-Witten theory,
with much more sophisticated foundations,
provides a simplifying completion of Hurwitz theory. 

\subsubsection{}

The present paper is the first of a series devoted to
the Gromov-Witten theory of target curves $X$. In 
subsequent papers, we will consider the equivariant
theory for $\proj^1$, 
the descendents of the
other cohomology classes of $X$,
and the connections to integrable hierarchies.
The equivariant Gromov-Witten theory of $\proj^1$ and 
the associated 2--Toda hierarchy will be the
subject of  \cite{OP}.

The introduction is organized as follows. 
We review the definitions of Gromov-Witten  and 
Hurwitz theory in Sections \ref{secgw} and \ref{hwz}. 
Shifted symmetric functions and
completed cycles are discussed in Section \ref{comcy}.
The basic GW/H correspondence is stated in Section \ref{gwhra}.

\subsection{Gromov-Witten theory}
\label{secgw}
The Gromov-Witten theory of a nonsingular target 
$X$ concerns integration
over the moduli space   $\overline{M}_{g,n}(X,d)$
     of stable degree $d$ 
maps from genus $g$, $n$-pointed curves to $X$. 
Two types of cohomology classes 
are integrated.
The {\em primary classes} are:
$$\text{ev}_i^*(\gamma) \in
H^2(\overline{M}_{g,n}(X,d), {\mathbb{Q}}),$$
where $\text{ev}_i$ is the morphism defined by evaluation at the 
$i$th marked point,
$$\text{ev}_i: \overline{M}_{g,n}(X)\rarr X,$$
and $\gamma\in H^*(X, {\mathbb{Q}})$.
The {\em descendent
classes} are:
$$ \psi_i^k \text{ev}_i^*(\gamma),$$
where
$\psi_i\in H^2(\overline{M}_{g,n}(X,d), {\mathbb{Q}})$
is
the first  Chern class of the cotangent line bundle $L_i$ on
the moduli space of maps.

Let $\omega \in H^2(X, {\mathbb{Q}})$ denote the Poincar\'e dual of the
point class. We will be interested here exclusively in the integrals of the
descendent classes of
$\omega$:
\begin{equation}
\label{llqq}
\lang \prod_{i=1}^n \tau_{k_i}(\omega) \rang_{g,d}^{\circ X} =
\int_{[\overline{M}_{g,n}(X,d)]^{vir}} \prod_{i=1}^n \psi_{i}^{k_i} \, 
\text{ev}_i^*(\omega).
\end{equation}
The theory is defined for all $d\geq 0$.

Let $g(X)$ denote the genus of the target.
The integral (\ref{llqq}) is defined to vanish unless the dimension constraint,
\begin{equation}
\label{wqqw}
2g-2+ d(2-2g(X))= \sum_{i=1}^n k_i,
\end{equation}
is satisfied.
If the subscript $g$ is omitted in the bracket notation $\lang \prod_i \tau_{k_i}(\omega)
\rang_d^X$,  the 
genus is specified
 by the dimension constraint from the remaining data. If the resulting genus is not an integer,
the integral is defined to vanish. Unless emphasis is required, the genus subscript
will be omitted.

The integrals (\ref{llqq}) constitute the {\em stationary sector} of
the
Gromov-Witten theory of $X$ since the images in $X$ of the
marked points are pinned by the
integrand. The total Gromov-Witten theory involves also  the
descendants of the identity and odd classes of $H^*(X, {\mathbb Q})$. 

The moduli space $\overline{M}_{g,n}(X,d)$ parameterizes
stable maps with connected domain curves.
However, Gromov-Witten theory may be also be defined with disconnected
domains.
If $C= \bigcup_{i=1}^l C_i$ is a disconnected curve with
connected components $C_i$, the arithmetic genus of $C$ is defined by:
$$g(C) = \sum_i g(C_i) - l+1,$$
where $g(C_i)$ is the arithmetic genus of $C_i$.
In the disconnected theory, the genus may be negative.
Let $\overline{M}^\bullet_{g,n}(X,d)$ denote the
moduli space of stable maps with possibly disconnected domains.

We will use the brackets brackets $\lang\ \rang^\circ$
as above in \eqref{llqq}
for integration in connected Gromov-Witten theory. 
The brackets
$\lang \ \rang^\bullet$ will be used for the
disconnected theory obtained by integration
against $[\overline{M}^\bullet_{g,d}(X,d)]^{vir}$.
The brackets
$\lang\ \rang$ will be used when it is not
necessary to distinguish between the connected and
disconnected theories.

\subsection{Hurwitz theory}
\label{hwz}

\subsubsection{}

The Hurwitz theory of a nonsingular curve $X$ concerns the
enumeration of covers of $X$ with specified ramification.
The ramifications are determined by the profile of the cover
over the branch points.

For Hurwitz theory, we will only consider covers,
$$\pi:C \rarr X,$$ where $C$ is nonsingular and $\pi$ is
dominant on each component of $C$. Let $d>0$ be the
degree of $\pi$.
The {\em profile} of $\pi$ over
a point $q\in X$ is
the partition $\eta$ of $d$ 
obtained from multiplicities of
$\pi^{-1}(q)$.

By definition,
a partition $\eta$  of $d$ is 
a sequence of integers,
$$
\eta = (\eta_1 \ge \eta_2 \ge \dots \ge 0),
$$
where  $|\eta|=\sum \eta_i = d$. Let
$\ell(\eta)$ denote the length of the partition $\eta$, 
and let $m_i(\eta)$ denote the multiplicity of the part $i$. 
The
profile of $\pi$ over $q$ is the partition $(1^d)$ if and
only if $\pi$ is unramified over $q$.

Let $d>0$, and 
let $\eta^1, \ldots, \eta^n$ be partitions of $d$ assigned to $n$ distinct
points
$q_1, \ldots, q_n$ of $X$.
A Hurwitz cover of $X$ of
 genus $g$, degree $d$, and monodromy $\eta^i$ at $q_i$ is
a morphism
\begin{equation}
\pi: C \rarr X\label{HurC}
\end{equation}
satisfying: 
\begin{enumerate}
\item[(i)] $C$ is a nonsingular curve of genus $g$,
\item[(ii)] $\pi$ has profile $\eta^i$ over $q_i$,
\item[(iii)]$\pi$ is unramified over 
$X \setminus
\{ q_1, \ldots, q_n\}$.
\end{enumerate}
Hurwitz covers may be considered with connected or
disconnected domains. 
The Riemann-Hurwitz formula,
\begin{equation}
\label{ggsspp}
2g(C)-2+ d(2-2g(X))  = \sum_{i=1}^n (d- \ell(\eta^i))\, ,
\end{equation}
 is valid for both connected
and disconnected Hurwitz covers. In the disconnected
theory, the domain genus may be negative. Since $g(C)$ is 
uniquely determined by the remaining data,  domain genus will be omitted
in the notation below.

Two covers 
$\pi: C \rarr X, \ \pi': C' \rarr X$ are isomorphic if
there exists an isomorphism of curves $\phi: C \rarr C'$ satisfying
$\pi'\circ \phi= \pi$. Up to isomorphism, there are only finitely many 
Hurwitz covers of $X$ of  genus $g$, degree $d$, and monodromy $\eta^i$ at $q_i$.
Each cover $\pi$ has a finite group of automorphisms
$\text{Aut}(\pi)$. 

The Hurwitz number, 
$$
\Hr_{d}^{X}(\eta^1, \ldots, \eta^n),
$$
is
defined to be the weighted count of the distinct,
\emph{possibly disconnected} Hurwitz covers $\pi$ with the
prescribed data. 
Each such cover is weighted by $1/|\text{Aut}(\pi)|$.

The GW/H correspondence is most naturally expressed as a relationship between the disconnected theories, hence 
the disconnected theories will be of primary interest to us.

\subsubsection{}

\label{xttx}
We will require an extended definition of Hurwitz numbers
valid in the degree $0$ case and in case the
ramification conditions $\eta$ satisfy $|\eta|\ne d$.  The 
Hurwitz numbers $\Hr_{d}^{X}$
are defined for all degrees $d\geq 0$
and all partitions $\eta^i$  by the following rules: 
\begin{enumerate}
\item[(i)] $\Hr_{0}^{X}(\emptyset, \ldots, \emptyset) = 1$,
where $\emptyset$ denotes the empty partition, 
\item[(ii)] if $|\eta^i| > d$ for some $i$ then the Hurwitz number vanishes, 
\item[(iii)] if $|\eta^i| \le  d$ for all $i$ then 
\begin{equation}
\Hr_d^{X}(\eta^1, \ldots, \eta^n) =
\prod_{i=1}^n \binom{m_1({\boldsymbol\eta}^i)}{m_1(\eta^i)} \cdot 
\Hr^{X}_{d}({{\boldsymbol\eta}}^1, \ldots, {{\boldsymbol\eta}}^n)\,.
\label{Hext}
\end{equation}
where ${\boldsymbol\eta}^i$ be the partition of size $d$
obtained from $\eta^i$ by adding $d-|\eta^i|$ parts of size 1.
\end{enumerate}

In other words, the monodromy condition $\eta$ at $q\in X$ with $|\eta|<d$
corresponds to counting  Hurwitz covers with
monodromy ${\boldsymbol\eta}$ at $q$ together with the data of 
a subdivisor of $\pi^{-1}(q)$ of profile $\eta$.

\subsubsection{}
The enumeration of Hurwitz covers of $\proj^1$
is classically known to be
equivalent to multiplication in the 
class algebra of the symmetric group. 
We review the theory here.

Let $S(d)$ be the symmetric group. Let 
$\Q S(d)$ be the group algebra. The {class algebra},
$$\ZZ (d) \subset \Q S(d),$$
is the center of the group algebra.

Hurwitz covers with profile 
$\eta^i$ over $q_i\in \proj^1$ 
canonically yield $n$-tuples of permutations
 $(s_1,\dots,s_n)$ defined {up to conjugation}
 satisfying:
\begin{enumerate}
\item[(i)] $s_i$ has cycle type $\eta_i$,
\item[(ii)] $s_1 s_2 \cdots s_n = 1.$
\end{enumerate}
The elements $s_i$ are determined by the
monodromies of $\pi$ around the points $q_i$.

Therefore, 
$\Hr_{d}^{\proj^1}(\eta^1, \ldots, \eta^n)$
equals the number of $n$-tuples satisfying conditions (i-ii) divided
by $|S(d)|$. The factor $|S(d)|$  accounts for over counting
and automorphisms.

Let $C_\eta\in\ZZ(d)$ be the conjugacy class corresponding
to $\eta$. We have shown:
\begin{align}
 \Hr_{d}^{\proj^1}(\eta^1, \ldots, \eta^n) &=
\frac1{d!} \, \left[C_{(1^d)}\right]\, \prod C_{\eta^i} \notag \\ &=
\frac1{(d!)^2} \tr_{\Q S(d)} \prod C_{\eta^i} \label{trS}
\end{align}
where $\left[C_{(1^d)}\right]$ stands for the coefficient of the identity
class and $\tr_{\Q S(d)}$ denotes the trace in the adjoint
representation. 

Let $\lambda$ be an
irreducible representation $\la$
of $S(d)$ of dimension $\dim \la$. 
The conjugacy class $C_\eta$ acts as a
scalar operator with eigenvalue
\begin{equation}\label{defff1}
\ff_\eta(\la)=|C_\eta| \, \frac{\chi^\la_\eta}{\dim\la} \,,
\quad |\la|=|\eta|\,,
\end{equation}
where $\chi^\la_\eta$ is the character of any element of 
$C_\eta$ in the representation $\la$.
The trace in equation  \eqref{trS} may be evaluated to yield the
basic character formula for Hurwitz numbers: 
\begin{equation}
  \label{Frob}
  \Hr_{d}^{\proj^1}(\eta^1, \ldots, \eta^n) =
\sum_{|\la|=d} \left(\frac{\dim\la}{d!}\right)^2 \, 
\prod_{i=1}^n \ff_{\eta^i}(\la) \,. 
\end{equation}

The character formula is easily generalized to include
the extended Hurwitz numbers (of Section \ref{xttx})
of
target curves $X$
of arbitrary genus $g$. 
The character formula can be traced to Burnside 
(exercise 7 in \S 238 of \cite{Burn}), see also\ \cite{D,J}.

Define $\ff_\eta(\la)$ for 
arbitrary partitions $\eta$ and irreducible
representations $\lambda$ of $S(d)$ by:
\begin{equation}\label{defff2}
\ff_\eta(\la)=\binom{|\la|}{|\eta|}\, 
|C_\eta| \, \frac{\chi^\la_\eta}{\dim\la} \,.
\end{equation}
If $\eta = \emptyset$, the formula is interpreted as:
$$
\ff_\emptyset(\la)=1 \,.
$$
For $|\eta|<|\la|$, the function $\chi^\la_\eta$ 
is defined
via the natural inclusion of symmetric groups $S(|\eta|) \subset S(d)$. 
If $|\eta|>|\la|$, the binomial in (\ref{defff2}) vanishes.

The character formula for extended Hurwitz numbers of genus $g$
targets $X$ is:
%
\begin{equation}
  \label{Burn}
  \Hr_{d}^{X}(\eta^1, \ldots, \eta^n) =
\sum_{|\la|=d} \left(\frac{\dim\la}{d!}\right)^{2-2g(X)} \, 
\prod_{i=1}^n \ff_{\eta^i}(\la) \,. 
\end{equation}

\subsection{Completed cycles}
\label{comcy}

\subsubsection{}

Let ${\mathcal P}(d)$ denote the set of partitions of $d$
indexing the irreducible representations of $S(d)$.
The Fourier transform, 
\begin{equation}
\label{qzzz}
\ZZ (d)\owns C_\mu \mapsto \ff_\mu\in 
\Q^{{{\mathcal P}(d)}}\,, \quad |\mu|=d\,,
\end{equation}
determines an isomorphism between
$\ZZ (d)$ and
the algebra of functions on ${\mathcal P}(d)$.
Formula \eqref{Frob} may be alternatively derived
as a consequence of the Fourier transform isomorphism.

Let ${\mathcal P}$ denote the set of all partitions
(including the empty partition $\emptyset$).
We may extend the Fourier transform 
(\ref{qzzz}) to define a map,
\begin{equation}\label{phi0}
\phi: \bigoplus_{d=0}^\infty \ZZ (d)\owns C_\mu \mapsto
\ff_\mu \in \Q^{\mathcal P} \,,
\end{equation}
via definition (\ref{defff2}).
%
The extended Fourier transform $\phi$ is no longer an isomorphism of 
algebras. However, $\phi$ is linear and injective.

We will see the image of $\phi$ in $\Q^{\mathcal P}$
is the algebra of \emph{shifted
symmetric functions} defined below
(see \cite{KO} and also \cite{OO}).

\subsubsection{}\label{sStab}

The shifted action of the symmetric group $S(n)$ on the algebra 
$\Q[\lam_1, \ldots, \lam_n]$ is defined by permutation of
the variables $\lam_i -i$.
Let $$\Q[\lam_1,\dots,\lam_n]^{*S(n)}$$ denote the invariants of the
shifted action. The algebra 
 $\Q[\lam_1,\dots,\lam_n]^{*S(n)}$
has a natural filtration by degree.

Define the algebra of shifted symmetric functions $\Ls$ in an 
infinite number of variables by 
\begin{equation}\label{defLs}
\Ls=\varprojlim \Q[\lam_1,\dots,\lam_n]^{*S(n)}\,,
\end{equation}
where the projective limit is taken 
in the category of filtered algebras 
with respect to the homomorphisms which send the last variable $\lam_n$ 
to $0$. 

Concretely, an element $\ff\in\Ls$ is
a sequence (usually presented as a series),
$$
\ff=\left\{\ff^{(n)}\right\}\,, \quad 
\ff^{(n)} \in \Q[\lam_1,\dots,\lam_n]^{*S(n)}\,,
$$
satisfying:
\begin{enumerate}
\item[(i)] the polynomials
$\ff^{(n)}$ are of  
uniformly bounded degree,
\item[(ii)] the polynomials $\ff^{(n)}$ are stable under restriction,
$$
\ff^{(n+1)}\big|_{\la_{n+1}=0} = \ff^{(n)}\,.
$$
\end{enumerate}
The elements of $\Ls$ 
will be denoted by boldface letters.

The  algebra $\Ls$ is filtered by degree.
The associated graded algebra $\gr\Ls$ is canonically isomorphic
to the usual algebra $\Lambda$ of symmetric functions as
defined, for example, in \cite{M}. 

A point $(x_1,x_2,x_3, \ldots)\in\Q^\infty$  is {\em finite} if all
but finitely many coordinates vanish. By construction, any
element $\ff\in\Ls$ has a well-defined evaluation at any
finite point. In particular, $\ff$ can be evaluated at any point
$$
\la=(\la_1,\la_2,\dots,0,0,\dots)\,,
$$
corresponding to a partition $\la$.  An elementary argument
shows functions
$\ff\in\Ls$ are uniquely determined by their values
$\ff(\la)$. Hence, $\Ls$ is canonically  a 
subalgebra of $\Q^{\mathcal P}$.

\subsubsection{}\label{spp}

The shifted symmetric power sum $\pp_k$ 
will play a central role in our study. Define
$\pp_k \in \Ls$ by:
\begin{equation}\label{defp}
\pp_k(\la)=\sum_{i=1}^\infty 
\left[(\la_i - i + \tfrac12)^k  - (- i + \tfrac12)^k\right]   + 
(1-2^{-k})\zeta(-k)\,. 
\end{equation}
The  shifted
symmetric polynomials,
$$
\sum_{i=1}^n 
\left[(\la_i - i + \tfrac12)^k  - (- i + \tfrac12)^k\right]
+ (1-2^{-k})\zeta(-k)
\,,
\quad n=1,2,3,\dots\,, 
$$
are of degree $k$ and are stable under restriction.
Hence, $\pp_k$ is well-defined.

The shifts by $\tfrac12$ in definition of $\pp_k$ appear
arbitrary --- their significance will be clear later. 
The peculiar $\zeta$-function constant term in $\pp_k$ will be
explained below.

The image of $\pp_k$ in 
$\gr \Ls \cong \Lambda$ is the usual $k$th power-sum
functions. Since the power-sums are well known to be 
free commutative generators of $\Lambda$, we conclude that
$$
\Ls=\Q[\pp_1,\pp_2,\pp_3,\dots] \,.
$$

The explanation of the constant term in \eqref{defp} is
the following. Ideally, we would like to define $\pp_k$
by 
\begin{equation} \label{ppxxx}
\pp_k \, \textup{``$=$''} \, \sum_{i=1}^\infty 
(\la_i - i + \tfrac12)^k\,.
\end{equation}
However, the above formula  violates stability and
diverges when evaluated at any partition
$\la$. In particular, evaluation at the empty partition $\emptyset$
yields:
\begin{equation}\label{ppkc}
\pp_k(\emptyset) \, \textup{``$=$''} \, \sum_{i=1}^\infty 
(- i + \tfrac12)^k\,. 
\end{equation}
Definition (\ref{ppxxx})
 can be repaired by subtracting the infinite constant \eqref{ppkc} inside
the sum in  
\eqref{defp} and  compensating by adding the 
$\zeta$-regularized value outside the sum.

The same 
regularization can be obtained in a more elementary
fashion by summing the following generating series:
$$
 \sum_{i=1}^\infty \sum_{k=0}^\infty  \frac{(-i+\frac{1}{2})^k z^k}{k!} =
\sum_{i=1}^\infty e^{z(-i+\frac{1}{2})}  = 
\frac{1}{z\,\cS(z)} \,,
$$
where, by definition, 
$$
\cS(z)= \frac{\sinh(z/2)}{z/2} = \sum_{k=0}^{\infty} 
\frac{z^{2k}}{2^{2k} \, (2k+1)!}\,. 
$$
The coefficients $c_i$ in the expansion, 
\begin{equation}
\frac{1}{\cS(z)} =
\sum_{i=0}^\infty c_i z^i\label{Bern}\, ,
\end{equation}
are essentially Bernoulli numbers. Since
$$
(1-2^{-k}) \, \zeta(-k) = k! \, c_{k+1} \,,
$$
the two above regularizations are equivalent.
The constants $c_{k}$ will play a important role.

It is convenient to
arrange the polynomials $\pp_k$ into a generating function:
\begin{equation}
\pp_k(\la) = k!\, [z^k]\, \ee(\la,z)\,,
\quad \ee(\la,z)= \sum_{i=0}^\infty e^{z(\la_i-i+\frac12)} \,,\label{ee}
\end{equation}
where $[z^k]$ denotes the coefficient of $z^k$ in the expansion 
of the meromorphic function $\ee(\la,z)$ in  Laurent series
about $z=0$.

\subsubsection{}

The function $\ff_\mu(\la)$, arising in the character formulas
for Hurwitz numbers, is shifted symmetric,
$$\ff_\mu \in \Ls,$$
a non-trivial result due to Kerov and Olshanski
 (see \cite{KO} and also 
\cite{OO,ORV}).  Moreover, the Fourier transform 
\eqref{phi0}
is a linear isomorphism,
\begin{equation}\label{phi}
\phi: \bigoplus_{d=0}^\infty \ZZ(d)\owns C_\mu \mapsto
\ff_\mu \in \Ls\,.
\end{equation}
The
identification of the highest degree term of
$\ff_\mu$ by
Vershik and Kerov (\cite{VK,KO}) yields:
\begin{equation}
\ff_\mu = \frac1{\prod \mu_i} \, \pp_\mu + \dots\,,
\label{fpl}
\end{equation}
where $\pp_\mu =\prod \pp_{\mu_i}$ and the dots stand
for terms of degree lower than $|\mu|$. 

The combinatorial interplay between the two 
mutually triangular linear
bases $\{\pp_\mu\}$ and $\{\ff_\mu\}$ of $\Ls$
is a fundamental aspect 
the algebra $\Ls$. In fact, these two bases will define the
GW/H correspondence.
 
Following \cite{EOZ}, we define the 
\emph{completed conjugacy classes} 
by
$$
\overline{C}_\mu=
\frac{1}{\prod_i \mu_i} \,
\phi^{-1}(\pp_\mu)  
\in \bigoplus_{d=0}^{|\mu|} \ZZ(d) \,.
$$
Since the basis $\{\pp_\mu\}$ is multiplicative, a
special role is played by the classes
$$
\overline{(k)} = \overline{C}_{(k)} \,, \quad k=1,2,\dots\,,
$$
which we call the \emph{completed cycles}. The formulas for first few
completed cycles are:
\begin{align*}
\overline{(1)}  = & (1) - \frac{1}{24} \cdot ()\, , \\ 
\overline{(2)} = & (2)\, , \\
\overline{(3)} = & (3) + (1,1)+ \frac{1}{12} \cdot (1) + \frac{7}{2880} 
\cdot ()\, ,\\
\overline{(4)} = & (4) + 2\cdot (2,1)+ \frac{5}{4} \cdot (2)\, ,
\end{align*}
where, for example, 
$$
(1,1) = C_{(1,1)} \in \ZZ (2) \,,
$$
is our shorthand notation for conjugacy classes.

Since $\ff_\mu(\emptyset)=0$ for any $\mu\ne\emptyset$,
the coefficient of the empty partition, 
$$()=C_\emptyset,$$ in $\overline{(k)}$ equals the constant
term of $\frac1k \pp_k$.

The 
{\em completion coefficients} $\rho_{k,\mu}$ determine the
expansions of the completed cycles,
\begin{equation}
\overline{(k)} = \sum_{\mu} \rho_{k,\mu} \cdot 
(\mu) \label{rho}\, .
\end{equation}
Formula \eqref{Bern} determining
the constants, 
$$
\rho_{k,\emptyset} = (k-1)! \, c_{k+1}  \,,
$$
admits a generalization determining
all the {completion coefficients},
\begin{equation}
  \rho_{k,\mu} = (k-1)!\frac{\prod \mu_i}{|\mu|!} \,
[z^{k+1-|\mu|-\ell(\mu)}] \, \cS(z)^{|\mu|-1} \,\prod \cS(\mu_i z) \,,
\label{rhof}
\end{equation}
where, as before, $[z^{i}]$ stands for the coefficient of 
$z^{i}$. Formula (\ref{rhof})
will be derived in Section \ref{srho}

The term {\em completed cycle} is appropriate as
$\overline{(k)}$ is obtained from $(k)$ by adding
non-negative multiples of conjugacy classes of strictly
smaller size (with the possible exception
of the constant term, which may be of either
sign). The non-negativity of $\rho_{k,\mu}$ for $\mu\neq \emptyset$
is clear from formula (\ref{rhof}). Also, the coefficient
$\rho_{k,\mu}$ vanishes unless the integer
$k+1-|\mu|-\ell(\mu)$ is even and non-negative.

We note 
the transposition 
$(2)$ is the unique cycle with no corrections required
for completion.

\subsubsection{}

The term {\em completed cycle} was suggested in \cite{EOZ}
when the
functions $\pp_k$ in \cite{BO,EO} were understood to
count degenerations
of Hurwitz coverings.
The GW/H correspondence explains the geometric meaning of
the completed cycles and, in particular, identifies the 
degenerate terms as contributions from the 
boundary of the moduli space of stable maps. 

In fact, completed cycles implicitly penetrate much 
of the theory of shifted symmetric functions. While the
algebra $\Ls$ has a very natural analog
of the Schur functions (namely, the \emph{shifted Schur
functions}, studied in \cite{OO} and many subsequent papers), 
there are several competing 
candidates for the analog of the power-sum symmetric
functions. The bases $\{\ff_\mu\}$ and $\{\pp_\mu\}$
are arguably the two finalists in this contest. The relationship
between these two linear bases can be studied using
various techniques, in particular, the methods of \cite{OO,ORV,LT}
can be applied. 


\subsection{The GW/H correspondence}
\label{gwhra}

\subsubsection{}

The
GW/H correspondence may be stated symbolically as:
\begin{equation}
  \label{taup}
\boxed{
  \tau_k(\omega) = \frac{1}{k!} \, \overline{(k+1)}}\,. 
\end{equation}
That is, descendents of $\omega$ are equivalent to completed
cycles.

Let $X$ be a nonsingular target curve.
The GW/H correspondence is the
following relation between the disconnected Gromov-Witten and 
disconnected Hurwitz theories: 
\begin{equation}
\label{bgwh}
\lang \prod_{i=1}^n \tau_{k_i}(\omega) \rang_{d}^{\bullet X} =
\frac{1}{\prod  k_i!} \, 
\Hr_d^{X}\left(\overline{(k_1+1)}, \ldots, \overline{(k_n+1)}
\right)\,,
\end{equation}
where the right-hand is defined by linearity via the
expansion of the completed cycles in ordinary conjugacy
classes.

The GW/H correspondence,
the completed 
cycle definition, and formula \eqref{Burn} together yield:
\begin{equation}
\label{Burn2}
\lang \prod_{i=1}^n \tau_{k_i}(\omega) \rang_{d}^{\bullet X} =
\sum_{|\la|=d} \left(\frac{\dim\la}{d!}\right)^{2-2g(X)} \, 
\prod_{i=1}^n \frac{\pp_{k_i+1}(\la)}{(k_i+1)!}\,.
\end{equation}
%
For $g(X)=0$ and $1$, the right side
can be expressed in the operator
formalism of the infinite wedge $\LV$ and explicitly evaluated, see 
Sections \ref{sP1} and \ref{sE}. 

The GW/H correspondence naturally extends
to relative Gromov-Witten theory, see Theorem \ref{fbgwh}.
In the relative context, the GW/H correspondence 
provides an invertible rule for exchanging descendent insertions $\tau_k(\omega)$
for
ramification conditions. 

The coefficients $\rho_{k,\mu}$  
are identified as connected $1$-point Gromov-Witten invariants
of $\proj^1$
relative to $0\in \proj^1$. 
The explicit formula \eqref{rhof} for the coefficients
is a particular case of the formula for $1$-point
connected GW invariants of $\proj^1$
relative to $0,\infty \in \proj^1$, see Theorem \ref{tm1pt}.

\subsubsection{}

Let us illustrate the GW/H correspondence in the special case of 
maps of degree $0$. In particular, we will see the role played
by the constants terms in the definition of $\pp_k$. 

In the degree $0$ case, the only partition 
$\la$ in the sum \eqref{Burn2} is the empty partition
$\la=\emptyset$. Since, by definition,
$$
\pp_k(\emptyset) = k! \, c_{k+1} \,,
$$
the formula \eqref{Burn2} yields
$$
\lang \prod \tau_{k_i}(\omega) \rang^{\bullet X}_{0} =
\prod c_{k_i+2}\,.
$$
The result is equivalent to the (geometrically obvious)
vanishing of all multipoint connected 
invariants, 
$$
\lang \tau_{k_1}(\omega) \cdots \tau_{k_n}(\omega)
\rang^{\circ X}_0 = 0 \,, \quad n>1\,,
$$
together with the following evaluation
of the connected degree 0, $1$-point function,
\begin{equation}
1+ \sum_{g=1}^\infty \lang \tau_{2g-2}(\omega) \rang^{\circ X}_{g,0}
\,  z^{2g} 
= \frac{1}{\cS(z)} \,.\label{d0gf}
\end{equation}
And, indeed, the result is correct, see \cite{FP,Pt}\,.

\subsubsection{}

A useful convention is to formally set the 
contribution $\lang \tau_{-2}(\omega) \rang^{\bullet X}_{0,0}$
of the unstable moduli space $\overline{M}_{0,1}(X,0)$ to 
equal $1$,
\begin{equation}
\lang \tau_{-2}(\omega) \rang^{\bullet X}_{0,0} = 1 \,.
\label{tau-20}
\end{equation}
This convention simplifies the form of the generating
function \eqref{d0gf} and several others functions in the paper.
In the disconnected
theory, the unstable contribution \eqref{tau-20} is allowed to
appear in any
degree and genus. 
Hence, in the disconnected theory, 
the convention is
equivalent to setting 
\begin{equation}
\tau_{-2}(\omega)  = 1 \,.
\label{tau-2}
\end{equation}
The parallel convention for the completed cycles
$$
\pp_0 = 0 \,, \quad \tfrac1{(-1)!}\,\pp_{-1} = 1 
$$
fits well with the formula \eqref{ee}\,.

\subsection{Plan of the paper}

\subsubsection{}

A geometric study of descendent integrals concluding with
a proof of the GW/H correspondence in the context of 
relative Gromov-Witten theory
is presented in Section \ref{motv}. 
The GW/H correspondence is Theorem \ref{fbgwh}.
A special case of GW/H correspondence
is assumed in the proof. The  
special case, the
GW/H
correspondence for the absolute Gromov-Witten theory of $\proj^1$, 
will 
be established by equivariant computations in \cite{OP}.

Relative Gromov-Witten theory is discussed in Section \ref{srel}. 
The completion coefficients \eqref{rho} are identified in Section \ref{ccyc}
as $1$-point Gromov-Witten invariants of $\proj^1$ relative to $0\in \proj^1$.

\subsubsection{}

The remainder of the paper deals with applications of the GW/H
correspondence. In particular, generating
functions for the stationary Gromov-Witten
invariants of targets of genus 0 and 1 
are evaluated.
These computations are most naturally executed in the
the infinite wedge formalism. We review the infinite representation
$\LV$
in Section \ref{aaa}. The formalism also provides
a convenient and powerful approach to the study of integrable 
hierarchies, see for example  \cite{K,MJD,SV}.

The stationary GW theory of $\proj^1$ relative to $0,\infty \in \proj^1$
is
considered in Section \ref{sP1}. We obtain a closed formula
for the corresponding $1$-point function in Theorem 
\ref{tm1pt}. The 
formula \eqref{rhof} for the completion coefficients is obtained
as a special case. A
generalization of Theorem \ref{tm1pt} for the $n$-point
function is given in Theorem \ref{tnpt}.

\subsubsection{}

The 2--Toda hierarchy
governing the Gromov-Witten theory of $\proj^1$ relative to 
$\{0,\infty\}\subset\proj^1$ is discussed in Section \ref{sToda}.
The main result is 
Theorem \ref{tmToda} which states the natural generating
function for relative GW-invariants is a $\tau$-function of 
the 2--Toda hierarchy of Ueno and Takasaki \cite{UT}.
Theorem \ref{tmToda} generalizes a 
result of \cite{Ot}. 

The flows of the Toda hierarchy are associated to the 
ramification conditions $\mu$ and $\nu$ imposed at $\{0,\infty\}$. 
The equations of the Toda 
hierarchy are equivalent to certain recurrence relations for
relative Gromov-Witten invariants, the simplest of which is made explicit in 
Proposition \ref{Tpr}.  

\subsubsection{}

The Gromov-Witten theory of $\proj^1$ was conjectured to be governed by
the Toda equation by  Eguchi and Yang \cite{EgY}, and 
also by Dubrovin \cite{Du}.
The Toda conjecture was 
further studied in in \cite{EgHY,EgHX,G,Ot,Pt}. 

The Toda conjecture naturally extends to
the $\C^\times$-equivariant Gromov-Witten theory of $\proj^1$.
We will prove in  \cite{OP} that
the equivariant theory of $\proj^1$ 
is governed by an integrable hierarchy which
can also be identified with the 2--Toda of \cite{UT}. 
The flows in the equivariant 2--Toda correspond to the insertions
of $\tau_k([0])$ and $\tau_k([\infty])$, where 
$$
[0],[\infty]\in H^*_{\C^\times}(\proj^1, Q)\,,
$$ 
are the classes of the torus fixed points. 

The
 equivariant 2--Toda hierarchy is \emph{different}
from the relative 2--Toda studied here.  However, the lowest equations
of both hierarchies agree on their common domain of applicability.

\subsubsection{}

In Section \ref{sE}, we discuss the stationary Gromov-Witten theory
of an elliptic curve $E$. The GW/H correspondence identifies the
$n$-point function of Gromov-Witten invariants of $E$
with the character of the infinite wedge representation
of $\gli$. This character has been previously computed in \cite{BO},
see also \cite{Oiw,EO}. We quote the results of \cite{BO} here and 
briefly discuss some of their implications, in particular, 
the appearance of quasimodular forms. 

While the GW/H correspondence is valid for 
all nonsingular target curves X, we do not know 
closed form evaluations for targets of genus $g(X)\geq 2$.
The targets $\proj^1$ and $E$ yield very beautiful theories.
Perhaps
the study of the Gromov-Witten theory of
higher genus targets will lead to the discovery of new structures. 

\subsection{Acknowledgments}
An important impulse for this work came from the results of \cite{EOZ}
and, more generally, from the line of research pursued in \cite{EO,EOZ}. 
Our interaction with S.~Bloch, A.~Eskin, and A.~Zorich played a very 
significant role in the development of the ideas presented here.

We thank E. Getzler and A. Givental for discussions of the 
Gromov-Witten theory of $\proj^1$, and T. Graber and Y. Ruan for discussions of the
relative theory.

A.O. was partially supported by 
DMS-0096246 and fellowships from the Sloan and Packard foundations.
R.P. was partially supported by DMS-0071473
and fellowships from the Sloan and Packard foundations.

%
%
%
%

\section{The geometry of descendents}
\label{motv}
\subsection{Motivation:  non-degenerate maps}
We begin by examining the relation between Gromov-Witten and Hurwitz
theory  in the context of non-degenerate maps with nonsingular domains.

Let $M^\bullet_{g,n}(X,d)\subset 
\overline{M}^\bullet_{g,n}(X,d)$
be the open locus of maps, $$\pi:(C,p_1, \ldots, p_n)\rarr
X,$$ where each connected component $C_i\subset C$ is
nonsingular and dominates $X$. 

Let $q_1,\ldots, q_n \in X$ be distinct points.
Define the closed substack $V$ by:
$$ V = \ev_{1}^{-1}(q_1) \cap \cdots \cap \ev_n^{-1}(q_n)
\subset M_{g,n}^\bullet(X,d)\,.
$$ 
The stacks $M^\bullet_{g,n}(X,d)$ and
$V$ are nonsingular Deligne-Mumford stacks of the  
expected dimensions --- see \cite{FanP}  for proofs.

The Hurwitz number $\Hr_{d}^{X}((k_1+1), \ldots,(k_n+1))$
may be defined by the 
enumeration of 
pointed Hurwitz covers 
 $$\pi:(C,p_1, \ldots, p_n)\rarr
(X,q_1,\dots,q_n)\,,$$
where
\begin{enumerate}
\item[(i)] $\pi(p_i)=q_i$,
\item[(ii)]
$\pi$ has ramification order $k_i$ at $p_i$.
\end{enumerate}
Here, $\pi$ has {\em ramification order} $k$ at $p$ if $\pi$
takes the local form $z \rarr z^{k+1}$ at $p_i$.
The count of pointed Hurwitz covers is 
weighted by $1/|\text{Aut}(\pi)|$ where $\text{Aut}(\pi)$ is the automorphism
group of the {\em pointed} cover.

The above enumeration of pointed covers
coincides with the definition of $\Hr_{d}^{ X}((k_1+1), \ldots, (k_n+1))$
given in Section \ref{hwz}.

\begin{pr}
\label{opres} Let $d>0$.
The algebraic cycle class, $$ 
\Big( \prod_{i=1}^n k_i! \, c_1(L_i)^{k_i}\ev_i^*(\omega)\Big) \   \cap [M^\bullet_{g,n}(X,d)]
\in A_0(
M^\bullet_{g,n}(X,d)),$$ is represented by the  
locus of
covers enumerated by
$\Hr_{d}^{X}((k_1+1), \ldots,(k_n+1))$. 
\end{pr}

\begin{proof}
Since $V$ represents $\prod_{i=1}^n \ev_i^*(\omega)$ in the Chow theory
of $M_{g,n}^\bullet(X,d)$, we may prove 
the locus of Hurwitz covers represents
$$\prod_{i=1}^n k_i!\, c_1(L_i)^{k_i}  \  \cap [V]$$
in the Chow theory of $V$.

First, consider the marked point $p_1$.
There exists a canonical section $s\in H^0(V, L_1)$
obtained from $\pi$ by the
following construction.
Let $\pi^*$ denote the pull-back map on functions:
\begin{equation}
\label{zzw}
\pi^*: m_{q_1}/m^2_{q_1} \rarr m_{p_1}/m^2_{p_1},
\end{equation}
where $m_{q_1}, m_{p_1}$ are the maximal ideals of
the points $q_1\in X$ and $p_1 \in C$ respectively.
Via the canonical isomorphisms,
$$m_{q_1}/m^2_{q_1} \stackrel{\sim}{=} T^*_{q_1}(X), \ \
 m_{p_1}/m^2_{p_1} \stackrel{\sim}{=} T^*_{p_1}(C),$$
the map (\ref{zzw}) is the dual of the 
differential
of $\pi$.
Since $q_1$ is fixed, the identification 
$m_{q_1}/m^2_{q_1} \stackrel{\sim}{=} \com$ yields
a section $s$ of $L_1$ by (\ref{zzw}).

The scheme theoretic zero locus $Z(s)  \subset V$ is easily seen
to be the (reduced) substack of maps where $p_1$ has ramification
order at least 1 over $q_1$. The cycle $Z(s)$ represents $c_1(L_1) \cap [V]$
in the Chow theory of $V$.

When restricted to
$Z(s)$, the pull-back of functions 
yields a map:
$$\pi^*: m_{q_1}/m^2_{q_1} \rarr m^2_{p_1}/m^3_{p_1}.$$
Hence, via the isomorphisms,
$$m_{q_1}/m^2_{q_1} \stackrel{\sim}{=} \com, \ \
 m^2_{p_1}/m^3_{p_1} \stackrel{\sim}{=} L_1^{\otimes 2},$$
a canonical section $s' \in H^0(Z(s), L_1^{\otimes 2})$ is
obtained. 
A direct scheme theoretic verification shows that $Z(s') \subset Z(s)$
is the (reduced) substack where $p_1$ has ramification order at least 
2 over $q$. Hence the cycle $Z(s')$ represents the cycle class
$2 c_1(L_1)^2$.

After iterating the
above construction, we find that $k_1! \, c_1(L_1)^{k_1}$ is
represented by the substack where $p_1$ has ramification order at least
$k_1$. At each stage, the reducedness of the zero locus is obtained
by a check in the versal deformation space of the ramified map 
(the issue of reducedness is local).

Since the cycles determined by ramification conditions at distinct markings $p_i$ are 
transverse, we conclude that
$\prod_{i=1}^n k_i! \, c_1(L_i)^{k_i} \cap [V]$ is represented by the
locus of Hurwitz
covers enumerated by
$\Hr_{d}^X((k_1+1), \ldots,(k_n+1))$.
\end{proof}

Proposition \ref{opres} shows a connection between descendent classes
and Hurwitz covers for the {\em open} moduli space $M_{g,n}^\bullet(X,d)$.
We therefore expect a geometric formula:
\begin{equation}
\label{mmww}
\lang
\tau_{k_1}(\omega)\cdots  \tau_{k_n}(\omega)\rang^{\bullet X}_{d} =
\frac{  \Hr_{d}^{X}((k_1+1), \ldots, (k_n+1)) }{
 \prod  k_i! }
 + \Delta,
\end{equation}
where $\Delta$ is a correction term obtained from the boundary,
$$\overline{M}^\bullet_{g,n}(X,d)\setminus M^\bullet_{g,n}(X,d).$$ 
The GW/H correspondence gives a description of this correction term 
$\Delta$.

For example, consider the case where $k_i=1$ for all $i$.
Then, since  2-cycles are already complete (see Section \ref{comcy}), 
the basic GW/H correspondence (\ref{bgwh}) yields an exact equality,
\begin{equation} 
\lang \tau_{1}(\omega)\cdots  \tau_{1}(\omega)\rang^{\bullet X}_{d} =
\Hr_{d}^{X}((2), \ldots, (2)),
\label{klk}
\end{equation}
which appears in \cite{Pt}.
However, the correction term 
$\Delta$ will not vanish in general.

We note that 
Proposition \ref{opres} holds for the connected moduli of maps and
connected Hurwitz numbers by the same proof. Since the disconnected
case will be more natural for the study of the correction equation 
(\ref{mmww}), the results have been stated in the disconnected case.

\subsection{Relative Gromov-Witten theory}\label{srel}

We will study the GW/H correspondence in the richer
context of the
Gromov-Witten
theory of $X$ relative to a finite set of distinct points $q_1,\ldots, q_m\in X$.
Let $\eta^1, \ldots, \eta^m$ be partitions of $d$. The moduli space
$$\overline{M}_{g,n}(X, \eta^1, \ldots, \eta^m)$$ parameterizes
genus $g$, $n$-pointed relative stable maps with monodromy
$\eta^i$ at $q_i$. Foundational
developments of relative Gromov-Witten theory in
symplectic and algebraic geometry can be found in
\cite{EGH,IP,LR,L}. 
The stationary sector of the relative 
Gromov-Witten theory is:
\begin{equation}
\label{neddde}
\lang  \prod_{i=1}^{n}
\tau_{k_i}(\omega),\eta^1, \ldots, \eta^m  \rang_{g,d}^{\circ X} =
\int_{[\overline{M}_{g,n}(X,\eta^1,\ldots,\eta^m)]^{vir}}
\prod_{i=1}^n \psi_{i}^{k_i} \,
\text{ev}_i^*(\omega),
\end{equation}
the integrals of descendents of $\omega$ relative to
$q_1, \ldots, q_m \in X$.

The genus \emph{and} the degree may be omitted 
in the notation (\ref{neddde}) as long as $m>0$. 
Again, the corresponding disconnected theory is denoted
by the brackets $\lang\ \rang^\bullet$. 

The stationary theory relative to $q_1, \ldots, q_m$ specializes
to the stationary theory relative to $q_1, \ldots, q_{m-1}$ when 
$\eta^m$ is the trivial partition $(1^d)$. In particular, when
all the partitions $\eta^i$ are trivial, the standard stationary
theory of $X$ is recovered.
A proof of this
specialization property is obtained from the degeneration formula
discussed in Section \ref{ddq} below.

The stationary Gromov-Witten theory of $\proj^1$ relative to
$0,\infty \in \proj^1$ 
will play a special role.
Let $\mu, \nu$ be partitions of $d$ prescribing the
profiles over $0, \infty \in \proj^1$ respectively.
We will use the notation,
\begin{equation}
\label{neddd}
\lang  \mu, \prod
\tau_{k_i}(\omega),\nu  \rang^{\proj^1},
\end{equation}
to denote integrals in the stationary theory of $\proj^1$
relative to $0,\infty \in \proj^1$.

\subsection{Degeneration} 
\label{ddq}
The degeneration formula for relative Gromov-Witten theory provides a formal
approach to the descendent integrals
$$\lang \tau_{k_1}(\omega)\cdots  \tau_{k_n}(\omega), \eta^1, \ldots,
\eta^m\rang^{\bullet X}_{d}.$$
Let $x_1, \ldots, x_n \in X$ be distinct fixed points.
Consider a family of curves with $n$ sections over the affine line,
$$\pi: ({\mathcal{X}},s_1, \ldots, s_n) \rarr {\mathbb A}^1,$$ 
defined by the following properties:
\begin{enumerate}
\item[(i)] $({\mathcal{X}}_{t}, s_1(t), \ldots, s_n(t))$ is 
isomorphic to the fixed
 data $(X, x_1, \ldots,x_n)$ for all $t\ne0$.
\item[(ii)] $({\mathcal{X}}_0, s_1(0), \ldots, s_n(0))$ is a comb consisting of
$n+1$ components (1 backbone isomorphic to $X$ and $n$ teeth  isomorphic to $\proj^1$). The
teeth are attached to the points 
$x_1,\ldots, x_n$ of the backbone. The section $s_i(0)$ lies on the $i$th
tooth.
\end{enumerate}
The degeneration $\pi$ can be easily constructed by blowing-up the $n$
points $(x_i,0)$ of the trivial family $X \times {\mathbb A}^1$.

The following result is obtained by viewing the family $\pi$ as 
a degeneration of the target 
in relative Gromov-Witten
theory.

\begin{pr}[\cite{EGH,IP,LR,L}] A 
degeneration formula holds for relative Gromov-Witten invariants:
\label{ffrrdd}
\begin{multline}\label{GWdgn}
  \lang \tau_{k_1}(\omega)\cdots \tau_{k_n}(\omega),\eta^1, \ldots,
  \eta^m
  \rang^{\bullet X}_d = \\
  \sum_{|\mu^1|, \ldots, |\mu^n|=d}\ \Hr^{X}_{d}(\mu^1, \ldots,
  \mu^n, \eta^1, \ldots, \eta^m) \prod_{i=1}^n \zz(\mu^i) \lang \mu^i,
  \tau_{k_i}(\omega) \rang^ {\bullet \proj^1}\,,
\end{multline}
where the sum is over all $n$-tuples $\mu^1,\dots,\mu^n$ of partitions
of $d$. 
\end{pr}
\noindent 

Here, the factor $\zz(\mu)$ is defined by: 
$$
\zz(\mu)= \left|\Aut(\mu)\right| \prod_{i=1}^{\ell(\mu)} \mu_i
$$
where $\Aut(\mu) \cong \prod_{i\geq 1} S({m_i(\mu)})$
is the symmetry group permuting equal parts of $\mu$. 
The factor $\zz(\mu)$ will occur often.

The right side of the degeneration formula \eqref{GWdgn}
involves the Hurwitz numbers and $1$-point stationary Gromov-Witten invariants
of $\proj^1$ relative to $0\in \proj^1$.
The degeneration formula together with the
definition of the Hurwitz numbers implies the specialization property of
relative Gromov-Witten invariants when $\eta^m= (1^d)$.

There exists an elementary analog of this degeneration 
formula in Hurwitz theory which yields:
\begin{multline}\label{Hdgn}
  \Hr_d^{X}
\left(\overline{(k_1)}, \ldots, \overline{(k_n)}, 
\eta^1, \ldots,\eta^m\right) = \\
\sum_{|\mu^1|, \ldots, |\mu^n|=d}  
\Hr^{X}_{d}\left(\mu^1, \ldots, \mu^n, \eta^1, \ldots, \eta^m\right)
\prod_{i=1}^n \zz(\mu^i) 
\, \Hr_d^{\proj^1}\left(\mu^i, \overline{(k_i)}\right)\,,
\end{multline}
where the sum is again over partitions $\mu^i$ of $d$. 

\subsection{The abstract GW/H correspondence}
\label{agwh}

Formula \eqref{GWdgn} can be restated as a substitution rule valid
in degree $d$:
\begin{equation} 
\tau_k(\omega)= \sum_{|\mu|=d} \left(\zz(\mu) \, \lang \mu,
  \tau_{k}(\omega) \rang^ {\bullet \proj^1}\right)\cdot (\mu) \,.
\label{abs1}
\end{equation}
The substitution rule replaces the descendents $\tau_k(\omega)$ by
ramification conditions in Hurwitz theory:
$$  \lang \tau_{k_1}(\omega)\cdots \tau_{k_n}(\omega),\eta^1, \ldots,
  \eta^m
  \rang^{\bullet X}_d = 
  \Hr^{X}_{d}(  
 -, \ldots,
  -,\eta^1, \ldots, \eta^m) \,.
$$
Hurwitz numbers on the right side are defined by inserting the respective
ramification
conditions
(\ref{abs1}) and expanding 
multilinearly. 
The substitution rule, however, is degree dependent by definition. 

A degree independent substitution rule is obtained by studying
the {\em connected} relative invariants. 
Disconnected invariants may be expressed as sums of products
of connected invariants obtained by all possible decompositions
of the domain and distributions of the integrand.
As the invariant  
$\lang \mu,
  \tau_{k}(\omega) \rang^ {\bullet \proj^1}$
a has single term in the integrand, an elementary argument yields:
%
\begin{equation}
\label{wwq}
\lang \mu, \tau_{k}(\omega) \rang^{\bullet \proj^1} =
\sum_{i=0}^{m_1(\mu)} \frac{1}{i!} 
\lang \mu - 1^i, \tau_{k}(\omega) \rang ^{\circ \proj^1},
\end{equation}
where $\mu-1^i$ denotes the partition $\mu$ with $i$ parts equal
to $1$ removed. 
Since
$$
\frac{\zz(\mu)}{i!} = \binom{m_1(\mu)}{i} \, \zz(\mu-1^i) \,,
$$
we may rewrite (\ref{wwq}) as:
\begin{equation*}
\zz(\mu)\lang \mu, \tau_{k}(\omega) \rang^{\bullet \proj^1} =
\sum_{i=0}^{m_1(\mu)} \binom{m_1(\mu)}{i} \, \zz(\mu-1^i)
\lang \mu - 1^i, \tau_{k}(\omega) \rang ^{\circ \proj^1}.
\end{equation*}
The following result is then obtained
from the definition of the extended Hurwitz numbers \eqref{Hext}.

\begin{pr} 
\label{qper}
A substitution rule 
for converting descendents 
to ramification conditions
holds: 
\begin{equation}
\tau_k(\omega)= \sum_{\nu} \left(\zz(\nu) \, \lang \nu,
  \tau_{k}(\omega) \rang^{\circ \proj^1}\right)\cdot (\nu) \,,
\label{abs2}
\end{equation}
where the summation is over {\em all} partitions $\nu$. 
\end{pr}

Proposition \ref{qper} is a degree independent,
abstract form of the GW/H correspondence. 
Clearly, only partitions $\nu$ of size at most 
$d$ contribute to the degree $d$
invariants. 
What remains
is the explicit identification of the coefficients in \eqref{abs2}.

\subsection{The leading term}
\label{acal}

Equating the dimension of the integrand in 
$\lang \nu, \tau_{k}(\omega) \rang^{\circ \proj^1}$ 
with the virtual dimension of the moduli space,
we obtain 
$$
k+1 = 2g-1 + |\nu| + \ell(\nu) \,.
$$
Since $g\ge 0$ and $\ell(\nu)\ge 1$, we find
$$
|\nu|\le k+1 \,.
$$
Moreover, $\nu=(k+1)$ is the only partition of size $k+1$
which actually appears in \eqref{abs2}. All other partitions
$\nu$ appearing in \eqref{abs2} have a strictly smaller
size. 

We will now determine the coefficient of $\nu=(k+1)$ in 
\eqref{abs2} by the method of  Proposition \ref{opres}.
The corresponding relative invariant is computed
in the following Lemma.

\begin{lm} \label{ggred}
For $d>0$, we have
$$
\lang (d), \tau_{d-1}(\omega) \rang^{ \proj^1}= \frac{1}{d!}\,.
$$
\end{lm}

\begin{proof} We first note the connected and disconnected invariants coincide,
$$\lang (d),\tau_{d-1}(\omega) \rang^{\circ \proj^1}=
\lang (d),\tau_{d-1}(\omega) \rang^{\bullet \proj^1}
,$$
since the imposed monodromy is transitive.
The genus of the domain is 0 by the dimension constraint.

Let $[\pi]\in \overline{M}_{0,1}(\proj^1, (d))$ be a stable
map relative to $0 \in \proj^1$,
$$\pi: (C, p_1) \rarr T \rarr \proj^1,$$
where $T$ is a destabilization of $\proj^1$ at $0$ and
$\pi(p_1) = \infty\in \proj^1$.
If $p_1$ lies on a $\pi$-contracted component $C_1 \subset C$
then, 
\begin{enumerate}
\item[(i)] $C_1$ must meet $\overline{ C \setminus C_1}$ in
at least 2 points by stability,
\item[(ii)] $\overline{ C \setminus C_1}$ must be connected
by the imposed monodromy at $0$.
\end{enumerate} Since conditions (i) and (ii) violate the
genus constraint $g(C)=0$, the marked
point $p_1$ is \emph{not}
allowed to 
lie on a $\pi$-contracted component of $C$.

The moduli space $\overline{M}_{0,1}(\proj^1, (d))$ is of expected
dimension $d$. 
By Proposition \ref{opres} pursued for relative maps, the cycle
$$ 
(d-1)! \, c_1(L_1)^{d-1} \ev^*_1(\omega) \cap [{M}_{0,1}(\proj^1,(d))] \in
 A_0( {M}_{0,1}(\proj^1,(d)))$$
is represented by the locus of covers enumerated by $\Hr_{0,d}((d),(d))$.

In fact, since $p_1$ does not lie on a $\pi$-contracted component of the domain
for any moduli point $[\pi]\in \ev_1^{-1}(\infty) 
\subset \overline{M}_{0,1}(\proj^1, (d))$, the
proof of Proposition \ref{opres} is 
is valid for the compact moduli space. The cycle 
$$ (d-1)! \, c_1(L_1)^{d-1} \ev^*_1(\omega)\cap 
[\overline{M}_{0,1}(\proj^1,(d))] 
 \in A_0(\overline
{M}_{0,1}(\proj^1,(d)))$$
is represented by the locus of covers enumerated by $\Hr_{0,d}((d),(d))$.

There is a unique cover $[\zeta]$ enumerated by $\Hr_{0,d}((d),(d))$. 
We may now complete the calculation:
\begin{eqnarray*}
\lang (d), \tau_{d-1}(\omega) \rang^{\bullet \proj^1} & = &
\int_ {[\overline{M}_{0,1}(\proj^1,(d))]}   c_1(L_1)^{d-1} \ev^*_1(\omega) \\
& = & \frac{1}{(d-1)!} \int_{[\zeta]} 1 \\ & = & \frac{1}{d!}, 
\end{eqnarray*}
since $[\zeta]$ is a cyclic Galois cover with automorphism group
of order $d$.
\end{proof}

Lemma \ref{ggred} provides an identification 
of the leading term in  the abstract GW/H correspondence \eqref{abs2}. 
\begin{cor}
  We have
  \begin{equation}
\tau_k(\omega) = \frac1{k!} \, {(k+1)} + \dots \,,
\label{leadt}
\end{equation}
where the dots stand for conjugacy classes $(\nu)$ with 
$|\nu|< k+1$. 
\end{cor}

\subsection{The full GW/H correspondence}
Let $X$ be a nonsingular curve.
The main result of the paper is a substitution rule for the relative
Gromov-Witten theory of $X$.

\begin{tm}\label{fbgwh}
A substitution rule
for converting descendents to ramification 
conditions holds: 
\begin{equation}\label{fGWH}
  \tau_k(\omega) = \frac1{k!} \, \overline{(k+1)} \,.
\end{equation}
\end{tm}

\noindent The full correspondence for the relative theory yields:
\begin{multline*}
 \lang \prod_{i=1}^n \tau_{k_i}(\omega), \eta^1, \ldots, \eta^m 
\rang_{d}^{\bullet X} = \\
\frac{1}{\prod  k_i!} \, 
\Hr_d^{X}\left(\overline{(k_1+1)}, \ldots, \overline{(k_n+1)}
,\eta^1, \ldots, \eta^m \right) \,.
\end{multline*}

Our proof of Theorem \ref{fbgwh} will
rely upon a special case  --- the case of the 
absolute Gromov-Witten theory of $\proj^1$. The formula, 
\begin{equation}
  \label{GWHp1}
 \lang \prod_{i=1}^n \tau_{k_i}(\omega) \rang_{d}^{\bullet \proj^1} =
\frac{1}{\prod  k_i!} \, 
\Hr_d^{\proj^1}\left(\overline{(k_1+1)}, \ldots, \overline{(k_n+1)}
\right) \,,
\end{equation}
will be proven in \cite{OP} as a result of equivariant computations.
We will now
deduce the general statement \eqref{fGWH} from \eqref{GWHp1}.

\begin{proof}
Let $\frac{1}{k!}\wtko$ denote the right side of the equality \eqref{abs2},
$$\frac{1}{k!}\wtko = 
\sum_{\nu} \left(\zz(\nu) \, \lang \nu,
  \tau_{k}(\omega) \rang^{\circ \proj^1}\right)\cdot (\nu) \,.$$
Define 
$\wtp_{k}$ by the 
the Fourier transform  \eqref{phi},
$$\phi\Big(   \wtk      \Big) = \frac{1}{k}\wtp_{k}.$$ 
 The equality \eqref{fGWH} is 
equivalent to the equality
\begin{equation}
\wtp_k \overset ? = \pp_k \,.\label{pfT1}
\end{equation}
As a result of \eqref{leadt}, we find:
\begin{equation}
\wtp_\mu = \pp_\mu + \dots  \,,\label{pfT2}
\end{equation}
where $\wtp_\mu = \prod \wtp_{\mu_i}$ and the dots stand for
lower degree terms. In other words, the transition matrix
between the bases $\{\wtp_\mu\}$ and $\{\pp_\mu\}$ is 
unitriangular. 

Let $\xi>0$ be a parameter, and let $l_\xi$ be the 
following linear form on the algebra $\Ls$:
$$
l_\xi(\ff) = \sum_{\la} \xi^{|\la|} \, \left(\frac{\dim \la}{|\la|!}
\right)^2 \, \ff(\la) \,.
$$
For example, $l_\xi(1)=e^\xi$. The series $l_\xi(\ff)$ is easily
seen define an entire function of $\xi$. 

The associated quadratic form,
\begin{equation}
(\ff,\gb)\mapsto l_\xi(\ff\cdot \gb) \,, \label{pfT3}
\end{equation}
is positive definite for any $\xi>0$.

Formula \eqref{Frob}, formula \eqref{GWHp1}, and the
definitions of the functions $\wtp_\mu, \pp_\mu$ yield
the equality, 
$$
l_\xi(\wtp_\mu) = l_\xi(\pp_\mu)\,,
$$
for all $\xi$ and $\mu$. In particular, we find
$$
l_\xi(\wtp_\mu \cdot \wtp_\nu) = l_\xi(\pp_\mu \cdot \pp_\nu)\,,
$$
for all $\mu$ and $\nu$. 
The transition 
matrix between the bases $\{\wtp_\mu\}$ and $\{\pp_\mu\}$
is therefore orthogonal with respect to the positive 
definite quadratic form \eqref{pfT3}.
By \eqref{pfT2}, the transition matrix is also unitriangular.
Hence, the transition is the identity and
equality is established in \eqref{pfT1}.
\end{proof}

\subsection{The completion coefficients}
\label{ccyc}
Theorem \ref{fbgwh} together with a 
comparison of the formulas \eqref{abs2} and \eqref{rho} yields
the following result.

\begin{pr}\label{tr1}
 The completion coefficients satisfy:
\begin{equation}
\frac{\rho_{k+1,\mu}}{k!} = \zz(\mu) \, \lang \mu,
  \tau_{k}(\omega) \rang^ {\circ \proj^1}\,.\label{rhoi}
  \end{equation}
In other words, the coefficients $\rho_{k,\mu}$ 
are determined by connected
relative $1$-point Gromov-Witten invariants of $\proj^1$ relative
$0\in \proj^1$.  
\end{pr}

The actual computation of these completion 
coefficients will be performed in Section \ref{sP1} using the operator 
formalism reviewed in Section \ref{aaa}. 
An explicit formula for the completion coefficients
will be given in Proposition \ref{rhofp}.

%
%
%
%

\section{The operator formalism}
\label{aaa}

The fermionic Fock space formalism reviewed here is a
convenient tool for manipulating the sums \eqref{Burn2}.
The operator calculus of the formalism is basic to the rest
of the paper. 
In Sections \ref{sP1} and \ref{sE},
the formalism is applied to 
the Gromov-Witten theory of targets of genus 0 and 1 respectively.
The formalism
underlies the study of the Toda hierarchy in Section \ref{sToda}.

\subsection{The infinite wedge}
\label{iwr}

\subsubsection{}

Let $V$ be a linear space with basis $\left\{\ul{k}\right\}$
indexed by the half-integers:
$$V = \bigoplus_{k \in \Z+ \sh} \com \, \ul{k}.$$
For each subset $S=\{s_1>s_2>s_3>\dots\}\subset \Z+\sh$ satisfying:
\begin{enumerate}
\item[(i)] $S_+ = S \setminus \left(\Z_{\le 0} - \sh\right)$ is
finite,
\item[(ii)]$S_- =
\left(\Z_{\le 0} - \sh\right) \setminus S$ is finite,
\end{enumerate}
we denote by  $v_S$ the following infinite wedge product:
\begin{equation}
v_S=\ul{s_1} \wedge \ul{s_2} \wedge  \ul{s_3} \wedge \dots\, .
\label{vS}
\end{equation}
By definition, 
 $$\LV= \bigoplus \com \,v_S
$$
is the linear space with basis $\{v_S\}$. 
Let $(\,\cdot\,,\,\cdot\, )$ be the inner product on  $\LV$ for which
$\{v_S\}$ is an orthonormal basis.

\subsubsection{}

The fermionic operator $\psi_k$ on $\LV$ is defined by wedge product with 
the vector $\ul{k}$, 
$$ 
\psi_k \cdot v = \ul{k} \wedge v  \,. 
$$
The operator $\psi_k^*$ is defined as the adjoint of $\psi_k$ 
with respect to the inner product  $(\,\cdot\,,\,\cdot\, )$. 

These operators satisfy the canonical anti-commutation relations:
\begin{gather}
  \psi_i \psi^*_j + \psi^*_i \psi_j = \delta_{ij}\,, \\
  \psi_i \psi_j + \psi_j \psi_1 = \psi_i^* \psi_j^* + \psi_j^*
  \psi_i^*=0.
\end{gather}
The {\em normally
ordered} products are defined by:
\begin{equation}\label{e112}
\nr{\psi_i\, \psi^*_j} =
\begin{cases}
\psi_i\, \psi^*_j\,, & j>0 \,,\\
-\psi^*_j\, \psi_i\,, & j<0 \,.
\end{cases}
\end{equation}

\subsubsection{} \label{sCH}

Let $E_{ij}$, for $i,j\in \Z+\sh$, be the standard basis of 
matrix units of $\gli$. 
The assignment 
$$
E_{ij} \mapsto\, \nr{\psi_i\, \psi^*_j}\ \ ,
$$
defines 
a projective representation of the Lie algebra $\gli=\mathfrak{gl}(V)$ on 
$\LV$. 

Normal ordering is introduced to avoid the infinite constants
which appear in the naive definition of the $\gli$-action 
on $\LV$. The ordering and divergence issues here are
closely related to the discussion in Section \ref{spp}.

For example, the 
action on $\LV$ of the identity matrix in $\gli$ is well-defined only
after normal ordering. 
Indeed, the operator,
$$
C= \sum_{k\in\Z+\frac12} \, E_{kk},
$$
corresponding to the identity matrix,
acts on the basis $v_S$ by:
$$
C \, v_S = (|S_+| - |S_-|) v_S \,.
$$
The operator $C$ is known as the \emph{charge} operator
\footnote{
The infinite wedge space is the mathematical formalization of
Dirac's idea of a sea of fermions filling all but finitely
many negative energy levels. The operator $C$ measures
the difference between the number $|S_+|$ of occupied 
positive energy levels (particles) and the number $|S_-|$ of
vacant negative energy levels (holes), whence the name.}.
The kernel of $C$, the zero charge subspace,
is spanned by the vectors 
\begin{equation*}
v_{\lam}=\ul{\lam_1-\tfrac12} \wedge \ul{\lam_2-\tfrac32} \wedge
\ul{\lam_3-\tfrac52} \wedge \dots 
\end{equation*}
indexed by all partitions $\la$. We will denote the kernel by $\LVc$. 

The operator  
$$
H= \sum_{k\in\Z+\frac12} k \, E_{kk}
$$
is called the \emph{energy operator}. The eigenvalues of $H$ 
on $\LVc$ are easily identified: 
$$ 
H \, v_\lambda = |\lambda| \, v_\lambda\,.
$$
The vacuum vector
$$
\vac = \ul{-\tfrac12} \wedge \ul{-\tfrac32} \wedge \ul{-\tfrac52}
\wedge \dots
$$
is the unique vector with the minimal (zero) eigenvalue of $H$.

\subsubsection{}

Define the translation operator $T$  by: 
\begin{equation}\label{defT}
T \, \ul{k_1} \wedge \ul{k_2} \wedge  \ul{k_3}  \wedge  \dots =
\ul{k_1+1} \wedge \ul{k_2+1} \wedge  \ul{k_3+1}  \wedge  \dots  \,.
\end{equation}
We see,
$$
T\,\psi_k \, T^{-1} = \psi_{k+1} \,, 
\quad  T\,\psi^*_k \, T^{-1} = \psi^*_{k+1}\,.
$$
We also find,
$$
T^{-1}\, C \, T  =  C + 1 \,.
$$
Hence, $T$ increases the charge by $1$.

\subsection{The operators $\cE$}
\label{scE}

\subsubsection{}

The operator $\cE_0(z)$ on $\LV$ is defined by:
\begin{equation}
\cE_0(z) = \sum_{k\in\Z+\frac12} \, e^{zk} \,  E_{kk} +  
\frac1{e^{z/2}-e^{-z/2}} \,,\label{E0}
\end{equation}
where the second term is a scalar operator on $\LV$. 
In fact,  
the scalar term in \eqref{E0} and the 
the constant term in \eqref{defp} have the same origin. 

Ideally, we would like $\cE_0(z)$ to be the naive
action on $\LV$ of
the following diagonal operator in $\gli$: 
$$
\ul{k} \mapsto e^{zk} \, \ul{k} \,.
$$
In other words, we would like to set 
\begin{equation}
\label{kkqqks}
\cE_0(z)  \, \textup{``$=$''} \, \sum_{k\in\Z+\frac12} \, e^{zk} \,  
\psi_k \, \psi_k^* \,,
\end{equation}
without normal ordering.
However, applied to the vacuum, definition (\ref{kkqqks}) yields:
\begin{equation}
\label{pqion}
\sum_{k=-\frac12,-\frac32,\dots} e^{zk} = \frac1{e^{z/2}-e^{-z/2}}\,,
\quad \Re z >0  \,,
\end{equation}
which may or may not make sense depending on $z$. 
We therefore define $\cE_0(z)$ using the normal ordering and 
then compensate by
adding the scalar \eqref{pqion} by hand.

In particular, we observe  
$$
\frac1{e^{z/2}-e^{-z/2}} =  \ee(\emptyset,z) \,,
$$
where the function $\ee(\la,z)$ is defined in \eqref{ee} and,
more generally, we find
\begin{equation}
  \label{Eeig}
  \cE_0(z) \, v_\la  = \ee(\la,z) \, v_\la \,.
\end{equation}
In other words, the functions $\ee(\la,z)$ are the eigenvalues
of the operator $\cE_0(z)$.

\subsubsection{}

Define the operators $\cP_k$ for $k>0$
by:
\begin{equation}
\cP_k = k! \, [z^k] \, \cE_0(z)\,,\label{defcP}
\end{equation}
where $[z^k]$ stands for the coefficient of $z^k$. From 
\eqref{ee} and \eqref{Eeig} we conclude:
\begin{equation}
  \label{Peig}
  \cP_k \, v_\la  = \pp_k(\la) \, v_\la \,.
\end{equation}
In particular, we find
$$
\cP_1 = H - \tfrac1{24}\,.
$$
The definition of the operators
$\cP_k$ is naturally extended as follows:
\begin{equation}\label{cP-1}
  \cP_{0} = C\,, \quad \tfrac1{(-1)!}\,\cP_{-1} = 1 \,.
\end{equation}
The extension is related to  convention \eqref{tau-2}.

\subsubsection{}

The translation operator $T$ acts on the operator
$\cE_0(z)$ by
\begin{equation}
  \label{TE}
 T^{-1} \, \cE_0(z) \, T = e^{z}  \,  \cE_0(z) \,.
\end{equation}
We find
\begin{equation}
  \label{TP}
 T^{-1} \, \frac{\cP_k}{k!}
 \, T = \sum_{m=0}^{k+1} \frac{1}{m!} \, \frac{\cP_{k-m}}{(k-m)!}\,,
\end{equation}
using  convention \eqref{cP-1}\,.

\subsubsection{}

For any $r\in\Z$, we define
\begin{equation}
\cE_r(z) = \sum_{k\in\Z+\frac12} \, e^{z(k-\frac{r}2)} \,  E_{k-r,k}
+ \frac{\delta_{r,0}}{\cs(z)} \, ,\label{Er}
\end{equation}
where the function $\cs(z)$ is defined by 
\begin{equation}
  \label{cs}
  \cs(z) = e^{z/2} - e^{-z/2} \,. 
\end{equation}
For $r\ne 0$, the normal ordering is not an issue 
and no constant term is required. The 
exponent in \eqref{Er} is set to satisfy: 
$$
\cE_r(z)^* = \cE_{-r}(z)^*\,,
$$
where the adjoint is with respect to the standard inner 
product on $\LV$. 

The operators $\cE$ satisfy the following fundamental 
commutation relation:
\begin{equation}\label{commcE}
  \left[\cE_a(z),\cE_b(w)\right] = 
\cs\left(\det  \left[
\begin{smallmatrix}
  a & z \\
b & w 
\end{smallmatrix}\right]\right)
\,
\cE_{a+b}(z+w)\,. 
\end{equation}
Equation \eqref{commcE} 
automatically incorporates the central 
extension of the $\gli$-action, which appears as the 
constant term in $\cE_0$ when $r=-s$.

\subsubsection{}

The operators $\cE$ specialize to the 
standard bosonic operators on $\LV$:
$$
\al_k = \cE_k(0)\,, \quad k\ne 0 \,.
$$
The commutation relation \eqref{commcE} specializes to 
the following equation
\begin{equation}
\label{qhhqq}
[\al_k, \cE_r(z)] =
 \cs(kz) \, \cE_{k+r}(z) \,.
\end{equation}
When $k+r=0$, equation (\ref{qhhqq}) has the following constant term: 
$$
\frac{\cs(kz)}{\cs(z)}=\frac{e^{kz/2}-e^{-kz/2}}{e^{z/2}-e^{-z/2}} \,.
$$
Letting $z\to 0$, we recover the standard relation:
\begin{equation}
[\al_k,\al_r]= k \, \delta_{k+r} \,. \label{commal}
\end{equation}

\subsubsection{}

The operators $\cE$ form a projective representation of 
the (completed) Lie algebra of differential operators on 
$\C^\times$, see for example \cite{KR,BO}. 

Let $x$ be the coordinate on $\C^\times$.
Identify $V$ with $x^{1/2} \C[x^{\pm1}]$ via the 
assignment
$$
\ul{k} \mapsto x^k \,.
$$
We then find the following correspondences:
$$
 \cP_k   \leftrightarrow  \left(x\frac{d}{dx}\right)^k    \,,
\quad \al_k \leftrightarrow x^k\,,
$$
where $x^k$ is considered as the operator of multiplication
by $x$. The correspondence is 
only a Lie algebra representation, and not a representation
of an associative algebra. 

The operator $\cE_0(z)$ corresponds to the following 
differential operator of infinite order
$$
T_z=\sum_{k} \frac{z^k}{k!} \, \left(x\frac{d}{dx}\right)^k\,,
$$
which acts on functions by rescaling their arguments:
$$
T_z\cdot f(x) =  f(e^z x) \,.
$$

\section{The Gromov-Witten theory of $\proj^1$}\label{sP1}

\subsection{The operator formula}

\subsubsection{}

The operator formalism will be used here to 
study the stationary Gromov-Witten invariants of
 of $\proj^1$ relative to $0,\infty\in\proj^1$,
\begin{equation*}
\lang  \mu, \prod 
\tau_{k_i}(\omega),\nu  \rang^{\bullet \proj^1},
\end{equation*}
and the corresponding connected invariants.
 
The GW/H correspondence \eqref{fGWH} together with \eqref{Frob} 
results in the following formula:
\begin{equation}\label{munu}
\lang  \mu, \prod_{i=1}^{n}
\tau_{k_i}(\omega),\nu  \rang^{\bullet \proj^1}=
\frac1{\zz(\mu) \zz(\nu)} 
\sum_{|\la|=|\mu|} \chi^\la_\mu \, \chi^\la_\nu \, \prod_{i=1}^{n}
\frac{\pp_{k_i+1}(\la)}{(k_i+1)!}  \,, 
\end{equation}
the derivation of which uses the equality
$$
|C_\mu| = |\mu|! / \zz(\mu) \,.
$$

\subsubsection{}

We first consider the following generating function,
\begin{equation}
  \label{Fd}
  F^\bullet_{\mu,\nu}(z_1,\dots,z_n) = 
\sum_{k_1,\dots,k_n=-2}^\infty 
\lang  \mu, \prod_{i=1}^{n}
\tau_{k_i}(\omega),\nu  \rang^{\bullet \proj^1} \,
\prod_{i=1}^n z_i^{k_i+1} \,,
\end{equation}
where convention \eqref{tau-2} is used for the $\tau_{-2}(\omega)$
insertions. Invariants in (\ref{Fd}) with $\tau_{-1}(\omega)$ insertions
are defined to vanish. 
Then, formula
\eqref{munu} may be rewritten as 
\begin{equation}
  \label{Fds}
  F^\bullet_{\mu,\nu}(z_1,\dots,z_n) =
\frac1{\zz(\mu) \zz(\nu)} 
\sum_{|\la|=|\mu|} \chi^\la_\mu \, \chi^\la_\nu \, \prod_{i=1}^{n}
\ee(\la,z_i)  \,.
\end{equation}
Our next goal is to recast the formula in terms of operators on $\LV$. 

\subsubsection{}

The following formula in $\LV$ is well-known:
$$
\prod_{i=1}^{\ell(\nu)}
\al_{-\nu_i} \, \vac = \sum_{|\la|=d} \, \chi^\la_\nu \,
v_\la \,, \quad |\nu|=d \, . 
$$
It is equivalent, for example, to the Murnaghan-Nakayama rule
for characters of symmetric group. Therefore,
using \eqref{Eeig}, we can express sum in the right side
of \eqref{Fds} as: 
$$
\left( \prod \cE_0(z_i) \, \prod \al_{-\nu_i} \, \vac, 
\prod \al_{-\mu_i} \, \vac\right) \,.
$$ 

For any operator $A$, we denote the diagonal matrix element of $A$ 
with respect to
the vacuum vector $\vac$ by angle brackets:
$$
\lang A \rang = (A \vac, \vac) \,.
$$
The above vacuum matrix element is the 
\emph{vacuum expectation}. Since, clearly, 
$$
\al_k^* = \al_{-k} \, ,
$$
formula \eqref{Fds} can be recast in the 
following operator form.

\begin{pr} We have 
\begin{equation}
  \label{Fdo}
  F^\bullet_{\mu,\nu}(z_1,\dots,z_n) =
\frac1{\zz(\mu) \zz(\nu)} 
\, \lang \prod_{i=1}^{\ell(\mu)} \al_{\mu_i}\, 
\prod_{i=1}^n \cE_0(z_i) \, 
\prod_{i=1}^{\ell(\nu)} \al_{-\nu_i}
\rang   \,.
\end{equation}
\end{pr}

\subsection{The 1-point series} \label{s1pt}

\subsubsection{}\label{s0pt}

We start by examining how the formalism
works in the (geometrically trivial) case of the 
$0$-point series. Formula \eqref{Fdo}
specializes to the following expression: 
\begin{equation}
  \label{0pt}
F^\bullet_{\mu,\nu}() =  \frac1{\zz(\mu) \zz(\nu)}
\, \lang \prod \al_{\mu_i}\, \prod \al_{-\nu_i}
\rang   \,.
\end{equation}
Observe, for positive $k$, the operator $\al_k$ annihilates the 
vacuum
$$
\al_k \, \vac = 0 \,, \quad k>0 \,.
$$
We can use the commutation relation \eqref{commal} 
repeatedly to move the operators $\al_{\mu_i}$ all 
the way to the right, after which the vacuum
expectation vanishes. Moving the operator $\al_{-\nu_i}$
all the way to the left has the same effect. 
Thus, a nonzero result is obtained
only in case all operators $\al_{\mu_i}$ and
$\al_{-\nu_i}$ annihilate in pairs via the 
commutation relation
\begin{equation}
[\al_k,\al_{-k}]= k \,. \label{commal2}
\end{equation}
This leads to the expected result 
$$
F^\bullet_{\mu,\nu}() = \frac{\delta_{\mu,\nu}}{\zz(\mu)} \,.
$$
{}From the geometric point of view, the 
commutation relation \eqref{commal2}, or the 
equivalent relation 
$$
\lang \al_k\,\al_{-k} \rang = k,
$$
is responsible for a $k$-fold covering of $\proj^1$ totally
ramified over $0$ and $\infty$. 

\subsubsection{}

Now we want to compute the $1$-point series
\begin{equation}
  \label{1pt}
F^\bullet_{\mu,\nu}(z) = \sum_{k=-2}^\infty 
\lang  \mu, 
\tau_{k}(\omega),\nu  \rang^{\bullet \proj^1} 
z^{k+1} \,,
\end{equation}
or, rather, the associated connected series
\begin{equation}
  \label{1ptc}
F^\circ_{\mu,\nu}(z) = \sum_{k=-2}^\infty 
\lang  \mu, 
\tau_{k}(\omega),\nu  \rang^{\circ \proj^1} \,
z^{k+1} \,. 
\end{equation}
We have
\begin{equation}
  \label{1pto}
F^\bullet_{\mu,\nu}(z) =  \frac1{\zz(\mu) \zz(\nu)}
\, \lang \prod \al_{\mu_i}\, \cE_0(z)\, \prod \al_{-\nu_i}
\rang   \,.
\end{equation}
We apply here the same strategy used 
to evaluate \eqref{0pt}: we move the operators
$\al_{\mu_i}$ to the right and move the operators 
$\al_{-\nu_i}$ to the left. 

We saw in the evaluation of the $0$-point series that the
commutation relation \eqref{commal2} accounts for 
a connected component without marked points.  
The commutators \eqref{commal2} make no contribution to
the connected series \eqref{1ptc}. 

All the action, therefore, happens as we commute the 
$\al$'s through the operator $\cE_0(z)$. The commutation 
is given by \eqref{qhhqq}. Applying formula \eqref{qhhqq} a total of
 $\ell(\mu)+\ell(\nu)$
times and using the obvious relation 
$$
|\mu| = |\nu|\,,
$$
we obtain the following result
\begin{align}
 F^\circ_{\mu,\nu}(z) & = 
\frac{\prod \cs(\mu_i z) 
\prod \cs(\nu_i z)}
{\zz(\mu) \, \zz(\nu)} \, \lang \cE_0(z) \rang \notag \\
&= \frac1{\zz(\mu) \zz(\nu)} \, 
\frac{\prod \cs(\mu_i z) 
\prod \cs(\nu_i z)}{\cs(z)} \label{1pt3}\,.
\end{align}
Using the function 
$$
\cS(z) = \frac{\cs(z)}{z} = \frac{\sinh z/2}{z/2}\,, 
$$
formula \eqref{1pt3} can be stated as follows. 
\begin{tm}\label{tm1pt}
For any two partitions $\mu$ and $\nu$ of the same size, we have 
  \begin{multline}
    \label{1ptf}
  \sum_{g=0}^\infty z^{2g}\, \lang  \mu, 
\tau_{2g-2+\ell(\mu)+\ell(\nu)}
(\omega),\nu  \rang^{\circ \proj^1} 
 =  \\ \frac1{|\Aut(\mu)|\,|\Aut(\nu)|} \, \frac{\prod \cS(\mu_i z) 
\prod \cS(\nu_i z)}{\cS(z)} \,. 
  \end{multline}
\end{tm}
Formula \eqref{d0gf} is recovered as the degree 0 case of Theorem
\ref{tm1pt}. More
generally, for
$\mu=\nu=(1^d)$, we obtain 
\begin{equation}
\sum_{g=0}^\infty z^{2g}\, \lang  
\tau_{2g-2+2d}
(\omega)  \rang^{\circ \proj^1} 
 =  \frac1{(d!)^2}\, \cS(z)^{2d-1}\,,
\label{prt1pt}
\end{equation}
which is the formula predicted in \cite{Pt} from the (then) conjectural
Toda equation. The Toda equation will be discussed in Section \ref{sToda}.
In particular, formula \eqref{prt1pt} can also be  
deduced from Proposition \ref{Tpr}.

\subsubsection{}

The product of $\cS$-functions
in \eqref{1ptf} satisfies an important property: the product
is symmetric in the combined set of variables
$\{\mu_i\}\cup\{\nu_i\}$.  This 
\emph{crossing symmetry} is very restrictive,
see \cite{EOZ}.
In particular, the symmetry implies that the full formula \eqref{1ptf} may
be obtained from the very special and 
degenerate case in which $\mu=(d)$. 
In fact, the property is
almost equivalent 
to the GW/H correspondence: the symmetry alone forces $\tau_k(\omega)$
to correspond to a \emph{linear} combination of the  $\pp_i$'s. 

We also observe that since $\cS(0)=1$, the coefficient of $z^{2g}$
in the product of $\cS$-functions in \eqref{1ptf} is well defined
as a symmetric functions of degree $2g$ in infinitely many variables.
In other words, we have the following stability: 
setting any variable to zero gives the analogous 
function in fewer variables, see the discussion in Section \ref{sStab}.

\subsubsection{}\label{srho}

{}From \eqref{1ptf} and Proposition \ref{tr1} we obtain 
the following result determining the completion coefficients.

\begin{pr}\label{rhofp} The completion 
coefficients \eqref{rho}  are given by 
\begin{equation}\label{rhoff}
  \rho_{k,\mu} = (k-1)! \frac{\prod \mu_i}{d!} \,
[z^{2g}] \, \cS(z)^{d-1} \,\prod \cS(\mu_i z) \,,
\end{equation}
where the $[z^{2g}]$ stands for the coefficient of 
$z^{2g}$ and the numbers $g$ and $d$ are defined by
$$
d=|\mu|\,, \quad k+1 = |\mu| +\ell(\mu) + 2g \,.
$$
\end{pr}

In particular, the terms for which $ |\mu| +\ell(\mu)$
reaches the maximal  value $k+1$ may be viewed together as  
principal terms of the completed cycle $\overline{(k)}$. 
For such terms, the genus $g$ vanishes and the
coefficient $\rho_{k,\mu}$ simply becomes
$$
\rho_{k,\mu} = (k-1)!\frac{\prod \mu_i}{|\mu|!}\,,
\quad |\mu| +\ell(\mu) = k+1 \,.
$$

The geometric interpretation of the coefficients
$\rho_{k,\mu}$ given in Proposition \ref{tr1} was \emph{not}
essential for the derivation of formula \eqref{rhoff}.

\subsection{The $n$-point series}

\subsubsection{}

The same strategy works for the evaluation of the 
general $n$-point series \eqref{Fd}, or, rather, the 
associated connected series 
$F^\circ_{\mu,\nu}(z_1,\dots,z_n)$. The result, however,
is somewhat more complicated to state. In particular, we 
require the following auxiliary function 
\begin{equation}
  \label{defG}
  G\left(
    \begin{matrix}
      a_1 & \dots & a_n \\
      z_1 & \dots & z_n 
    \end{matrix}\right) = 
\lang \cE_{a_1}(z_1)\, \dots \, \cE_{a_n}(z_n) \rang^\circ \,,
\end{equation}
where the superscripted circle indicates the connected part
of the vacuum expectation, that is, 
$$
\lang \cE_{a_1}(z_1)\, \cE_{a_2}(z_2) \rang^\circ =
\lang \cE_{a_1}(z_1)\, \cE_{a_2}(z_2) \rang -
\lang \cE_{a_1}(z_1)\rang \, \lang \cE_{a_2}(z_2) \rang\,,
$$
et cetera. The function \eqref{defG} 
clear vanishes unless the condition
$$
a_1 + \cdots + a_n = 0 
$$
is satisfied. Also, the equation
$$
G\left(
    \begin{matrix}
      0 \\
      z 
    \end{matrix}\right) = \frac1{\cs(z)}  \, 
$$
is clear.

\subsubsection{}

For $n>1$, the function \eqref{defG} can be computed 
recursively as follows. First, if $a_1\le 0$ then 
\eqref{defG} vanishes
$$
  G\left(
    \begin{matrix}
      a_1 & \dots & a_n \\
      z_1 & \dots & z_n 
    \end{matrix}\right)=0\,, \quad a_1 \le 0 \,.
$$
If $a_1>0$, then by commuting the operator $\cE_{a_1}(z_1)$
all the way to the right using the commutation relation 
\eqref{commcE}, we obtain
\begin{multline*}
  G\left(
    \begin{matrix}
      a_1 & \dots & a_n \\
      z_1 & \dots & z_n 
    \end{matrix}\right) =  \\
\sum_{i=2}^n 
\cs\left(\det  \left[
\begin{matrix}
  a_1 & a_i \\
z_1 & z_i
 \end{matrix}\right]\right)
\, 
G\left(
    \begin{matrix}
      a_2 & \dots & a_i+a_1 & \dots & a_n \\
      z_2 & \dots & z_i+z_1 & \dots & z_n 
    \end{matrix}\right)\,, \quad a_1>0 \,.
\end{multline*}
The above rules can be easily converted into a non-recursive form. 
For example, for $n=2$ we have 
$$
  G\left(
    \begin{matrix}
      a & -a \\
      z_1 & z_2 
    \end{matrix}\right) = 
\begin{cases}
    {\displaystyle\frac{\cs(a(z_1 + z_2))}{\cs(z_1+z_2)}} \,, & a>0\,, \\
0 & a \le 0 \,. 
\end{cases}
$$

\subsubsection{}

The strategy of Section \ref{s1pt} applies to 
evaluation of \eqref{Fdo} with minor modification.  
Each of the $\al$'s now has a {\em choice}
of operator $\cE$ with which to interact (that is, with which to 
commute). This choice can be conveniently formalized
in terms of a function 
$$
f:\{\mu_i\}\cup \{-\nu_i\} \to \{1,\dots,n\} \,,
$$
where $\{\mu_i\}\cup \{-\nu_i\}$ is 
considered as a multiset, 
that is, a set with possible repetitions. The evaluation $f(\mu_i)=j$ indicates
the commutator of $\al_{\mu_i}$
with $\cE_0(z_j)$ is taken. 

%
\begin{tm}\label{tnpt}
Let $M$ denote the multiset $\{\mu_i\}\cup \{-\nu_i\}$.
 We have
  \begin{multline}\label{nptf}
F^\circ_{\mu,\nu}(z_1,\dots,z_n) = \\
\frac1{\zz(\mu) \zz(\nu)} \,
\sum_f \left(\prod_{m\in M} \cs\left(|m|\, z_{f(m)}\right)
\right) \,
G\left(
    \begin{matrix}
      &\dots & \sum_{f^{-1}(i)} m & \dots&  \\
      &\dots & z_i & \dots &  
    \end{matrix}\right)\,,
\end{multline}  
where the sum is over all functions $f: M \to \{1,\dots, n\}$. 
\end{tm}

Since the summation over $f$ in \eqref{nptf} involves 
$n^{\ell(\mu)+\ell(\nu)}$
terms, formula \eqref{nptf} is only effective
if the partitions $\mu$ and $\nu$ have few parts. For 
partitions of large length, especially
for the case $\mu=\nu=(1^d)$, a more effective answer is
given by the Toda equations, see in particular 
Proposition \ref{Tpr}.

%
%
%
%

\section{The Toda equation}\label{sToda} 

We study here the Toda equations for 
the relative Gromov-Witten theory of $\proj^1$. The Toda equations
are equivalent to certain recurrence relations 
for the relative invariants.

\subsection{The $\tau$-function}
\label{tfun}

The $\tau$-function is 
a generating function of the relative 
invariants of $\proj^1$ which is convenient from the point
of view of integrable hierarchies. 

\subsubsection{}

Let $t_1,t_2,\dots$ be a sequence of indeterminates.
Consider the following \emph{vertex operators}:
$$
\Gamma_\pm(t) = \exp\left(\sum_{k>0} t_k \, \frac{\al_{\pm k}}{k} 
\right)\, .
$$
We easily obtain:
$$
\Gamma_-(t) = \Gamma_+(t)^*, 
$$
$$
\Gamma_+(t) = \sum_{\mu} \frac{t_\mu}{\zz(\mu)} \, \prod \al_{\mu_i} \,,
$$
where $t_\mu = \prod t_{\mu_i}$. 
The above sum is over all 
partitions $\mu$.

\subsubsection{}

Define the $\tau$-function 
for Gromov-Witten theory of $\proj^1$ relative to $0,\infty\in\proj^1$ by:
\begin{equation}
\tau_{\proj^1}(x,t,s) = 
\sum_{|\mu|=|\nu|} t_\mu \, s_\nu \,
\lang \mu, \exp\left(\sum_{i=0}^\infty x_i\, \tau_i(\omega)\right), \nu
\rang^{\bullet\proj^1} \,,\label{tauP1}
\end{equation}
where $x_0,x_1,\dots$ is new set of variables.
The following conventions will hold for the degree 0 constant terms:
$$
\lang \emptyset, \emptyset \rang^\bullet_0 = 1 \,,
\quad \lang \emptyset, \emptyset \rang^\circ_0 = 0 \,.
$$
The $\tau$-function is often called the {\em partition function}.

Formula
\eqref{Fdo} and the definition \eqref{defcP} of the 
operators $\cP_k$ together yield an operator formula for
$\tau_{\proj^1}$.

\begin{pr}\label{ptau1} We have
$$
\tau_{\proj^1}(x,t,s) = \lang \Gamma_+(t) \, 
\exp\left(\sum_{k=0}^\infty \frac{x_k}{(k+1)!} \, \cP_{k+1} \right) \,
\Gamma_-(s) \rang \,.
$$  
\end{pr}

\subsubsection{}

By the usual relation between the connected and 
disconnected theories, the logarithm of $\tau$ generates the 
connected invariants: 
\begin{align}
  \cF_{\proj^1}(x,t,s) &= 
 \sum_{|\mu|=|\nu|} t_\mu \, s_\nu \,
  \lang \mu, \exp\left(\sum_{i=0}^\infty x_i\, \tau_i(\omega)\right),
  \nu \rang^{\circ\proj^1}
\notag \\ &= \ln \tau_{\proj^1}(x,t,s)\label{FP1} \,,
\end{align}
where the first equality is the definition of the function 
$\mathcal{F}_{\proj^1}$. The function $\mathcal{F}_{\proj^1}$
is known as the \emph{free energy}.

\subsection{The string equation}

\subsubsection{}

As a slight extension of stationary Gromov-Witten theory, 
we allow the appearance of 
$$
\tau_0(1) \,,
$$
a marked point with 
no imposed conditions. The $\tau_0(1)$-insertions
are known as \emph{punctures}. In both the connected and disconnected theory,
the insertions of $\tau_0(1)$ can be removed using {\em string equation}:
\begin{equation}
\lang \tau_0(1) \prod_i \tau_{k_i}(\omega) \rang = 
\sum_j \lang \prod_i \tau_{k_i-\delta_{i,j}}(\omega) \rang
\label{string}\,.
\end{equation}
By the same principle, one removes any number of punctures,
which can be expressed as follows: 
\begin{equation}\label{estring}
\lang e^{y \tau_0(1)} \prod_i \tau_{k_i}(\omega) \rang = 
\lang \prod_i \left(\sum_{m\ge 0} \frac{y^m}{m!}\,
\tau_{k_i-m}(\omega) \right)\rang \,.
\end{equation}

\subsubsection{}

In the standard interpretation of the string equation, all 
the negative descendants are set to zero.
Also, there is the following unique exception to the string 
equation in the connected theory:
$$
\lang \tau_0(1)^2 \, \tau_0(\omega) \rang_{0,0} = 1 \,.
$$
In the disconnected theory, of course, the exception 
propagates in all degrees and genera. 

An equivalent way of managing the exceptional case
is to declare the string equation always valid,
while simultaneously changing the interpretation of 
the output. Recall our conventions for the
disconnected stationary theory: 
\begin{equation}
\tau_k(\omega) =
\begin{cases}
  1 \,, & k = -2 \,,\\
  0 \,, & k \ne 2,\, k<0 \,,
\end{cases}
\label{tauneg}
\end{equation}
We now observe the following interpretation 
of the string equation is equivalent to the 
standard one:
\begin{enumerate}
\item[(i)] we first apply the string 
equation, without exceptions and 
without setting $\tau_{-1}(\omega)$ and $\tau_{-2}(\omega)$
to zero, 
repeatedly to 
remove all $\tau_0(1)$-insertions,
\item[(ii)] after which we apply the rules \eqref{tauneg} to the 
resulting stationary Gromov-Witten invariant\,.
\end{enumerate}

\subsubsection{}

The form of the string equation is unchanged in 
relative Gromov-Witten theory. 
Let us add an additional string variable $y_0$ to the 
generating function \eqref{tauP1}:
\begin{multline}
\tau_{\proj^1}(x,t,s,y_0) = \\
\sum_{|\mu|=|\nu|} t_\mu \, s_\nu \,
\lang \mu, \exp\left(y_0\,\tau_0(1)+
\sum_{i=0}^\infty x_i\, \tau_i(\omega)\right), \nu
\rang^{\bullet\proj^1} \,.\label{tauP2}
\end{multline} 
Similarly, the function 
$$
\cF_{\proj^1}(x,t,s,y_0) = \ln \tau_{\proj^1}(x,t,s,y_0)
$$
is the generating function for connected invariants
in the presence of punctures.

Equations \eqref{estring} and \eqref{TP} together
with the commutation of operator $T$ with 
vertex operators $\Gamma_\pm$ results in the 
following generalization of Proposition \ref{ptau1}:

\begin{pr}\label{ptau2} For $n\in\Z$, we have
$$
\tau_{\proj^1}(x,t,s,n) = \lang T^{-n} \, \Gamma_+(t) \, 
\exp\left(\sum_{k=0}^\infty \frac{x_k}{(k+1)!} \, \cP_{k+1} \right) \,
\Gamma_-(s) \, T^n\rang \,.
$$  
\end{pr}

\subsection{The Toda hierarchy}

\subsubsection{}

By a standard argument (which can be found, for 
example, in \cite{Ot} and will be explained in 
more detail in \cite{OP}), Proposition \ref{ptau2}
yields the following result.

\begin{tm}\label{tmToda} The sequence 
$$
\{\tau_{\proj^1}(x,t,s,n)\}\,, \quad n\in \Z \,,
$$
is a $\tau$-function of the 2--Toda hierarchy of Ueno and 
Takasaki {\rm \cite{UT}} in the variables $t$ and $s$. 
In particular, the lowest equation of 
this hierarchy is:
\begin{equation}\label{lt}
  \frac{\partial^2}{\partial t_1 \partial s_1} \log \tau(n) =
  \frac{ \tau(n+1) \,\tau(n-1)} { \tau(n)^2} \,,
\end{equation}
where $\tau(n)=\tau_{\proj^1}(x,t,s,n)$\,. 
\end{tm}

The two sequences of flows in this hierarchy are 
connected with two 
ramification conditions $\mu$ and $\nu$ in the relative
Gromov-Witten theory, and \emph{not} with the descendent 
insertions $\tau_k(\omega)$.

In particular,  since 
$$
\tau_1(\omega)=\overline{(2)}=(2)\,,
$$
the function $\tau_{\proj^1}$
specializes under the restriction
 $$x_2=x_3=\dots=0$$ 
to the $\tau$-function of \cite{Ot} enumerating Hurwitz covers
with arbitrary branching over $0,\infty\in\proj^1$ and simple
ramifications elsewhere. Thus, Theorem \ref{tmToda}
generalizes the results of \cite{Ot}.

\subsubsection{}\label{altT}

A 2--Toda hierarchy of a different kind arises in 
in the equivariant GW theory
of $\proj^1$, see \cite{OP}. The flows of the equivariant
2--Toda hierarchy are associated to the 
insertions of $\tau_k([0])$ and $\tau_k([\infty])$, where 
$$
[0],[\infty]\in H^*_{\C^\times}(\proj^1)
$$ 
are the classes of the $\C^\times$-fixed 
points in the equivariant cohomology of $\proj^1$. 

In 
the non-equivariant limit, both $[0]$ and $[\infty]$ yield
the point class $\omega$. In the non-equivariant specialization,
the 2--Toda becomes a 1--Toda hierarchy for the absolute stationary 
 Gromov-Witten theory of $\proj^1$ described by 
the function 
$$
\tau^{\textup{abs}}_{\proj^1}(x,q)=
\tau_{\proj^1}(x,t_1,0,0,\dots,s_1,0,0,\dots) \,,
$$
where the variable $q=t_1 s_1$ keeps track of degree. 
On this absolute stationary submanifold, the lowest
equations of the two
different Toda hierarchies coincide, see Section \ref{sttT}.

Getzler in \cite{G} has
constructed an extension of the $1$--Toda hierarchy, 
the extra flows of which correspond to
the descendents
$\tau_k(1)$. 
In other words, Getzler's hierarchy  describes the full 
absolute non-equivariant Gromov-Witten theory of $\proj^1$.
Getzler has proven the extended hierarchy is essentially
 equivalent to the union of the stationary $1$--Toda hierarchy 
and the Virasoro constraints \cite{G}.

\subsubsection{}

In terms of the free energy  \eqref{FP1}, equation \eqref{lt}
reads
\begin{equation}\label{lt2}
  \frac{\partial^2}{\partial t_1 \partial s_1} \cF_{\proj^1}(x,t,s,y_0) =
\exp\left(
\Delta \, \cF_{\proj^1}(x,t,s,y_0)\right)
\,,
\end{equation}
where $\Delta$ is the following divided difference operator in the 
string variable $y_0$
$$
\Delta f(y_0) = f(y_0+1) - 2 f(y_0) + f(y_0-1) \,.
$$
Equivalently, the operator $\Delta$ can be interpreted as the 
insertion of 
$$
e^{\tau_0(1)} - 2 + e^{-\tau_0(1)} = \cs\left(\tau_0(1)\right)^2 \,.
$$
The exponential on the right side of \eqref{lt2} can be 
interpreted as a generating function for disconnected invariants,
modified by the action of the operator $\Delta$. 

Observe that, by the definition of $\cF_{\proj^1}$, the 
coefficient of $t_\mu s_\nu$ in the expansion of 
$\frac{\partial^2}{\partial t_1 \partial s_1} \cF_{\proj^1}(x,t,s)$
is equal to:
$$
(m_1(\mu)+1)(m_1(\nu)+1)\, \lang 
\mu+1, \exp\left(\sum_{i=0}^\infty x_i\, \tau_i(\omega)\right),
  \nu+1 \rang^{\circ\proj^1}\,,
$$
where $\mu+1$ denotes the partition $\mu\cup\{1\}$.

\subsubsection{}

By the string equation \eqref{estring}, 
the effect of the operator $\Delta$ on 
an $n$-point function is the following:
\begin{multline}\label{Del-npt}
  \sum_{k_i} \lang \mu, \cs\left(\tau_0(1)\right)^2 \prod
\tau_{k_i}(\omega) , \nu \rang \, \prod z_i^{k_i+1} = \\
%
%
\cs\left(\sum z_i \right)^2 
\sum_{k_i} \lang \mu, \prod 
\tau_{k_i}(\omega) , \nu \rang \, \prod z_i^{k_i+1} \,. 
\end{multline}
In particular,  the result vanishes when $n=0$. Hence,  the 
$0$-point functions do not appear in the right-hand side
of \eqref{lt2}. 

We may now translate equation \eqref{lt2} to the following relation
for $n$-point functions.

\begin{pr}\label{Tpr}
The Toda equations \eqref{lt}, \eqref{lt2}
are equivalent to the following recurrence relation for 
$n$-point functions. For any $\mu$ and $\nu$ of the same size, we have
\begin{multline}
 F^\circ_{\mu+1,\nu+1}
(z_1,\dots,z_n) = \\ 
\frac1{(m_1(\mu)+1)(m_1(\nu)+1)}\sum_{\{(S_i,\mu^i,\nu^i)\}}
\prod_i  \cs\left(\Sigma_{S_i}\right)^2  \, 
F^\circ_{\mu^i,\nu^i}\left(z_{S_i}\right)\,,
\end{multline}
where the summation is over all sets of triples 
$$
\{(S_i,\mu^i,\nu^i)\}\,, 
$$
such that $\{S_i\}$ is a partition of the set $\{1,\dots,n\}$
into nonempty disjoint subsets:
$$
\{1,\dots,n\} = \bigsqcup S_i\,, \quad S_i\ne\emptyset\,, 
$$
similarly, $\{\mu^i\}$ and $\{\nu^i\}$ satisfy
$$
\mu = \bigcup \mu^i\,, \quad \nu = \bigcup \nu^i \,,\quad
|\mu_i|=|\nu_i|\,,
$$
and where, by definition,
 $z_S=\{z_i\}_{i\in S}$ and  $\Sigma_S=\sum_{i\in S} z_i$.
\end{pr}

It is instructive to notice the consistency of this result 
with the result of Theorem \ref{tm1pt}. 

\subsubsection{}\label{sttT}

We will now consider the absolute stationary Gromov-Witten theory of $\proj^1$.
From the generating function $\tau_{\proj^1}$, the absolute
 specialization 
$\tau_{\proj^1}^{\textup{abs}}$
is obtained by setting
$$
t_2=t_3=\dots=s_2=s_3=\dots=0 \,.
$$
The restricted function $\tau_{\proj^1}^{\textup{abs}}$ depends on $t_1$ and $s_1$
only through the weight $(t_1 s_1)^d$ multiplying terms of 
degree $d$. Similarly, its dependence on the variable $x_0$ is
exclusively through the weight $e^{x_0(d-\frac1{24})}$ multiplying 
terms of degree $d$ in $\tau_{\proj^1}^{\textup{abs}}$. The constant term 
$-\frac1{24}$ can be transformed into an overall factor of 
$e^{-x_0/24}$. Since
$$
\frac{\partial^2}{\partial x_0^2} \log e^{-x_0/24} = 0  \,,
$$
we see
$$
t_1 s_1\,\frac{\partial^2}{\partial t_1 \partial s_1 } \, \log
\tau_{\proj^1}^{\textup{abs}} =
\frac{\partial^2}{\partial x_0^2} \log
\tau_{\proj^1}^{\textup{abs}} \,.
$$
We now replace the $\frac{\partial^2}{\partial t_1 \partial s_1 }$
derivative in \eqref{lt2} by derivatives with respect with $x_0$. Then, we
we set $t_1s_1=q$. We obtain the following result.

\begin{pr} The generating function 
$$
\cF^{\textup{abs}}_{\proj^1}(x,y_0,q) = 
 \sum_d q^d
  \lang \exp\left(y_0\,\tau_0(1)+\sum_{i=0}^\infty x_i\, \tau_i(\omega)\right)
  \rang^{\circ\proj^1}_d 
$$
for the absolute invariants of $\proj^1$ satisfies
the following version of the 
Toda equation \eqref{lt2} 
\begin{equation}\label{lt3}
  \frac{\partial^2}{\partial x_0^2}\, \cF^{\textup{abs}}_{\proj^1}(x,y_0,q) =
q \exp\left(
\Delta \, \cF^{\textup{abs}}_{\proj^1}(x,y_0,q)\right)
\,.
\end{equation}
\end{pr}

In contrast to \eqref{lt}, \eqref{lt2}, the 
differentiation in \eqref{lt3} is with respect to the
variable coupled to the insertion of $\tau_0(\omega)$. 
Equation \eqref{lt3} is the lowest equation in another
Toda hierarchy, mentioned in Section \ref{altT},
the flows of which are associated to  
descendent insertions in the absolute Gromov-Witten theory
of $\proj^1$.

\section{The Gromov-Witten theory of $E$}\label{sE}

Since a nonsingular cubic can be degenerated to a nodal rational curve, 
the degeneration principle explained in Section \ref{ddq} yields
the following expression for 
the Gromov-Witten invariants of an elliptic curve $E$ 
in terms of relative invariants of $\proj^1$:
\begin{equation}
\lang \prod \tau_{k_i}(\omega)
\rang^{\bullet E}_{d} = \sum_{|\mu|=d} \zz(\mu) \,
\lang \mu, \prod \tau_{k_i}(\omega), \mu
\rang^{\bullet \proj^1}\,.\label{degE}
\end{equation}
Here the sum is taken over all partitions $\mu$ of $d$. 

Consider the following $n$-point generating function 
$$
F_E(z_1,\dots,z_n;q) = \sum_{d\ge 0} q^d \sum_{k_1,\dots,k_n} 
\lang  \prod_{i=1}^{n}
\tau_{k_i}(\omega)  \rang_{d}^{\bullet E} \,
\prod_{i=1}^n z_i^{k_i+1} \,,
$$
which includes contributions of all degrees. From the 
degeneration formula \eqref{degE} and the 
operator formula \eqref{Fdo}, we conclude that
\begin{align}
 F_E(z_1,\dots,z_n;q) &= \sum_{\mu} \frac{q^{|\mu|}}{\zz(\mu)} \,
 \, \lang \prod \al_{\mu_i}\, \prod \cE_0(z_i) \, \prod \al_{-\mu_i}
\rang \notag \\
&= \tr_0 \, q^H \prod \cE_0(z_i) \label{trE}\,,
\end{align}
where $\tr_0$ denotes the trace in the charge zero
subspace $\LVc\subset\LV$, spanned by the vectors
$v_\la$ or, equivalently, by the vectors
$$
\prod \al_{-\mu_i} \, \vac \,,
$$
as $\la$ or $\mu$ range over all partitions. 
The 
vectors $\prod \al_{-\mu_i} \, \vac$ are orthogonal with 
norm squared equal to $\zz(\mu)$, see Section \ref{s0pt}. 
Also,  the energy operator $H$ in \eqref{trE} was defined in 
Section \ref{sCH}. 

The trace \eqref{trE} has been previously computed in 
\cite{BO}, see also \cite{Oiw,EO}. The result is the
following. Introduce the product
$$
(q)_\infty = \prod_{n=1}^\infty (1-q^n) \,,
$$
and the genus 1 theta function
$$
\Th(z)=\Th_{\frac12,\frac12}(z;q)= 
\sum_{n\in\Z} (-1)^n q^{\frac{(n+\frac12)^2}{2}} e^{(n+\frac12)z} \,.
$$ 
Up to normalization, $\Th(z)$ is the only odd genus 1 theta function ---
the normalization is immaterial 
as the formula  will be homogeneous in $\Th$. 

\begin{tm}[\cite{BO}]
We have
\begin{multline}\label{npoint}
F_E(z_1,\dots,z_n;q) = \\ \frac1{(q)_\infty}
\sum_{\substack{\textup{all $n!$ permutations} 
\\\textup{of $z_1,\dots,z_n$}}}
\frac
{\det\left[
\dfrac{\Th^{(j-i+1)}(z_1+\dots+z_{n-j})}{(j-i+1)!}
\right]_{i,j=1}^n}
{\Th(z_1) \, \Th(z_1+z_2) \cdots \Th(z_1+\dots+z_n)}\,,
\end{multline}
where in the $n!$ summands the $z_i$'s are permuted in all
possible ways.
\end{tm}

\noindent Here, $\Th^{(k)}$ denotes the $k$-th derivative of $\Th$.
If $k<0$, the standard convention $1/k!=0$ is followed. Hence, negative 
negative derivatives do not appear in formula \eqref{npoint}.

A qualitative conclusion which may be drawn 
is that the $z$-coefficients of 
\eqref{npoint} are quasimodular forms in the degree variable 
$q$. Concretely, for any collection of the $k_i$'s, we have
\begin{equation}
(q)_\infty \, \sum_{d=0}^\infty q^d \, 
\lang  \prod
\tau_{k_i}(\omega)  \rang_{d}^{\bullet E} \in 
\Q[E_2,E_4,E_6]_{\sum (k_i+2)} \,,
\label{qmod}
\end{equation}
where $\Q[E_2,E_4,E_6]$ denotes the ring (freely) 
generated by the Eisenstein series
$$
E_k(q)=\frac{\zeta(1-k)}2+\sum_n  
\left(\sum_{d|n} d^{k-1}\right) q^n
$$
of weight $k=2,4,6$, and the lower index specifies
the homogeneous component of weight $\sum (k_i+2)$. 
This quasimodularity condition is both 
very useful and very restrictive. The modular
transformation  
relates the $q\to 1$ behavior of the series \eqref{qmod}
with its $q\to 0$ behavior, thus connecting large degree
invariants with low degree invariants. 

Since the $2$-cycle is complete,
$$
\tau_1(\omega) = \overline{(2)} = (2)
$$
the quasimodularity (\ref{qmod}) generalizes the quasimodularity
of generating functions for simply branched coverings
of the torus studied in \cite{D,KZ}. 

Further discussion of the properties of the function 
\eqref{npoint} can be found in \cite{BO,EO}. In particular,
\cite{EO} contains the asymptotic analysis of 
this function as $q\to 1$, which corresponds to 
the $d\to\infty$ asymptotics of the GW-invariants.

\vspace{+10 pt}
\noindent
Department of Mathematics \\
UC Berkeley \\
Berkeley, CA 94720\\
okounkov@math.berkeley.edu \\

\vspace{+10 pt}
\noindent
Department of Mathematics \hfill Department of Mathematics\\
\noindent California Institute of Technology \hfill  Princeton University\\
\noindent Pasadena, CA 91125 \hfill Princeton, NJ 08544\\
\noindent rahulp@cco.caltech.edu \hfill rahulp@math.princeton.edu
\end{document}